\documentclass[12pt,fleqn]{article}
\usepackage{comment}
\usepackage{graphicx}
\usepackage{amsmath}
\usepackage{hyperref}
\usepackage{float}
\usepackage{color}
\usepackage{xcolor}
\usepackage{tikz}
\usepackage{subcaption}
\usetikzlibrary{arrows}
\begin{document}

\title{Effect of diagnostic testing on the isolation rate in a compartmental model with asymptomatic groups
	\thanks{The work of Zuzana Chladn\'a was partly supported by Slovak Grant Agency APVV-0096-12, and the work of Jana Kopfov\'a was supported by the institutional support for the development of research organizations I\v C 47813059.}
}
%



\author{Zuzana Chladn\'a         \and
        Jana Kopfov\'a          \and
        Dmitrii Rachinskii      \and
        Pavel \v{S}tep\'anek
        etc.
}
\date{Received: date / Accepted: date}
\maketitle

\begin{abstract}
    We present two epidemiological models, which extend the classical SEIR model
    by accounting for the effect of
    indiscriminate quarantining, isolation of infected individuals based on testing and the presence of asymptomatic individuals. 
    Given a constraint that limits the maximal number of simultaneous active cases,
we demonstrate that the isolation rate, which enforces this constraint, 
decreases with the increasing testing rate.
The models predict that massive testing allows to control the infection spread using a much lower isolation rate than in the case of indiscriminate quarantining.
We also show that dynamics of infection are more sensitive to the parameters of the asymptomatic groups
than the corresponding parameters of the symptomatic groups.
\end{abstract}


\section{Introduction}
As the race for developing effective COVID-19 vaccines and the search for medications for the treatment of the disease continues, home quarantine measures remain the main means for controlling the epidemic.
Restrictions imposed by home quarantine policies on businesses and individuals  help to decrease the number of social contacts and 
thus slow down the spread of the disease.
These restrictions differed in their scale, forms and timings in the states and countries where they were implemented by the government and health authorities. 
However, massive lockdowns are not a sustainable solution for economic, social and psychological reasons. The situation is complicated by the asymptomatic COVID-19 infection. The asymptomatic individuals can account for as many as
40-45\% of infections \cite{topol} and  can transmit the virus to others for an extended period of time.

Asymptomatic infection, the latent infectious phase of the disease and cases with mild symptoms 
can be detected by diagnostic tests. 
As such, massive diagnostic testing of the population and isolation of positively tested individuals can potentially offer an alternative to 
indiscriminate quarantining of large population groups 
\cite{mina}. 
The isolation strategy based on diagnostic screening was implemented by a number of organizations, business and government authorities.
For example, major airlines such as American Airlines and Luft Hansa added on-site pre-flight rapid testing facilities at a number of American and German airports.
Duke University comprehensive COVID-19 testing program received results from 16,146 tests administered to students and faculty from November 7\,--\,13, 2020.  The Chinese city of Qingdao reported testing its entire population of nine million people for COVID-19 over a period of five days in October, 2020.  Slovakia tested 3.6 million people -- two thirds of its population -- in
two days on October 31 -- November 1, 2020 and repeated 
testing the same population in a week's time in the first attempt of a large-scale blanket testing campaign in Europe.
These massive testing efforts were assisted by new pulling technologies and the development of affordable rapid tests 
such as Abbott's BinaxNOW.


The health and government authorities are concerned with keeping the number of active cases below the level dictated by the  capacity of the 
health care system.
Several recent studies have attempted to predict and analyze the effect of quarantine measures on dynamics of COVID-19 pandemic using compartmental models of mathematical epidemiology.
A variant of the standard well mixed SEIR model 
\begin{equation}
\label{SEIR}
	\begin{array}{l}
	\dot S= -\beta SI,\\
	\dot E= \beta SI - \omega E,\\
    \dot I=\omega E - \delta I,\\
	 \dot R=\delta I,\\
	\end{array}
\end{equation}
involving additionally the transmission from exposed to susceptible individuals,
was used to analyze the  effect of social distancing and reducing the number of contacts \cite{SIRQ2}. 
Another implementation of the quarantine policy 
was studied in \cite{SIRQ1},  where the authors adapted 
the SIR model 
assuming that all the infected individuals are isolated  
after the incubation period. Stochastic   age-structured   transmission   models were applied  to   explore  a    range   of   intervention 
 scenarios \cite{kucharski}.   A model with 
a two-threshold switching prevention strategy 
 predicted that flattening the curve can lead to periodic recurrence of the disease \cite{our}. The authors of 
 all  the above mentioned studies 
 conclude that 
 the interventions are effective in reducing the infection peak; however, 
extreme interventions  are  likely to be  required  to  contain the infection spread. 

Amid the ongoing efforts of
bringing the COVID-19 epidemic under control continues, it is important to explore different scenarios using a variety of models and modeling assumptions.
In this paper, we are interested in the potential of massive diagnostic testing measures for reducing the quarantining rate, which enforces a given constraint on the maximum of active cases.
A few 
earlier
studies addressed the role of testing in the trajectory of the epidemic. In particular,
the discrete time adaptation of the SEIR model proposed in \cite{testing} suggests that the isolation based on  testing interventions
can make unnecessary costly lockdown measures.

Below we propose two compartmental models, which account for indiscriminate quarantine measures, 
targeted detection and isolation of infected individuals based on diagnostic testing and a possible difference of the epidemiological parameters of the symptomatic and asymptomatic groups.
We compute the trajectories of the infected and quarantined populations
in the case when indiscriminate quarantining is applied and in the case when targeted quarantining
is facilitated by diagnostic testing. We then evaluate the relative efficiency of the testing effort 
by matching the parameters of both scenarios to ensure that they produce the same infection peak.
As the measure of the efficiency, the time-integral of the total quarantined population over the duration of the epidemic is used.
Finally, we discuss the sensitivity and robustness of trajectories to the epidemiological parameters 
of the symptomatic and asymptomatic groups.

The models are presented in the next section.
For simplicity, we do not consider the possibility of reinfection, even though a few cases of COVID-19 reinfection have been seemingly confirmed. 
Contact tracing per se is not included in the models 
but can be accounted for as a factor increasing the isolation rates. Two different types of tests are used to detect the presence of the virus  
during the illness and the presence of antibodies after the illness. 
We consider the former type only. 
The results are presented in Section \ref{sec3}.
The values  of the epidemiological parameters used in the numerical simulations are based on the average characteristics of COVID-19 that we distilled from the current literature.

\textcolor{red}{\bf 
%
%
}

\section{Models}\label{sec2}

\subsection{Basic model} 
\begin{figure}
\begin{center}
\begin{tikzpicture}
[
  font=\sffamily,
  every matrix/.style={ampersand replacement=\&,column sep=1.5cm,row sep=2cm},
  source/.style={draw,thick,rounded corners,fill=yellow!20,inner sep=.25cm,text width=1.8cm},
  to/.style={->,>=stealth',shorten >=1pt,semithick},
  every node/.style={align=center}]
 
  \matrix{
    \node[source](S){$S$};\&\&; \node[source](SQ){$S_Q$}; \\
    \node[source](E){$E$};\&\&; \node[source](EQ){$E_Q$}; \\
    \node[source](I_a){$I_a$};
    \&;\node[fill=none,draw=none](belowI_a){v};\node[source](I_{sQ}){$I_{sQ}$}; \&\node[fill=none,draw=none](belowLQ){};\node[source](LQ){$I_{aQ}$};\\
     \node[source](R){$R$};\&\& \node[source](RQ){$R_Q$};\&\\
  };
  \draw[to] (S) -- node[midway,left] {$(\beta$,$\beta(1-\rho))$} (E);
    \draw[to] (E) -- node[midway,left] {$k\omega$} (I_a);
     \draw[to] (SQ) -- node[midway,left] {$(\beta(1-\rho),\beta(1-\rho)^2) $} (EQ);
     \draw[to] (S) -- node[midway,above] {$\chi$} (SQ);
      \draw[to] (E) -- node[midway,above] {$\chi$} (EQ);
   \draw[to] (E) -- node[midway,left] {\hspace{-0.4cm} $(1-k)\omega$} (I_{sQ});
      \draw[to] (EQ) -- node[midway,left] { $(1-k)\omega$ \hspace{0.1cm}} (I_{sQ});
  \draw[to] (EQ) -- node[midway,left] {$k\omega$} (LQ);
  \draw[to] (I_a) -- node[midway,left] {$\delta$} (R);
  \draw[to] (R) -- node[midway,above] {$\chi$} (RQ);
    \draw[to] (I_a) -- node[midway,above] {$\psi$} (I_{sQ});
  \draw[to] (LQ) -- node[midway,left] {$\delta$} (RQ);
  \draw[to] (LQ) -- node[midway,above] {$\psi$} (I_{sQ});
  \draw[to] (I_{sQ}) -- node[midway,below] {\hspace{-0.2cm} $\delta$} 
 (R);
 \path[->]
 (I_a) edge [bend right=45] node[midway,above] {$\chi$} (LQ);

 \draw [loosely dashed,ultra thick] (-2.3,5.4) rectangle (4.75, -6);
 
  
\node at (4.3,-5.5) {\large $Q$};
\end{tikzpicture}

\end{center}
   \caption{\label{Diagram1} Compartmental diagram of model \eqref{model1}. Arrows show the flows (with the corresponding rate  parameters) between the  compartments. The dashed rectangle includes the quarantined populations. }  
\end{figure}
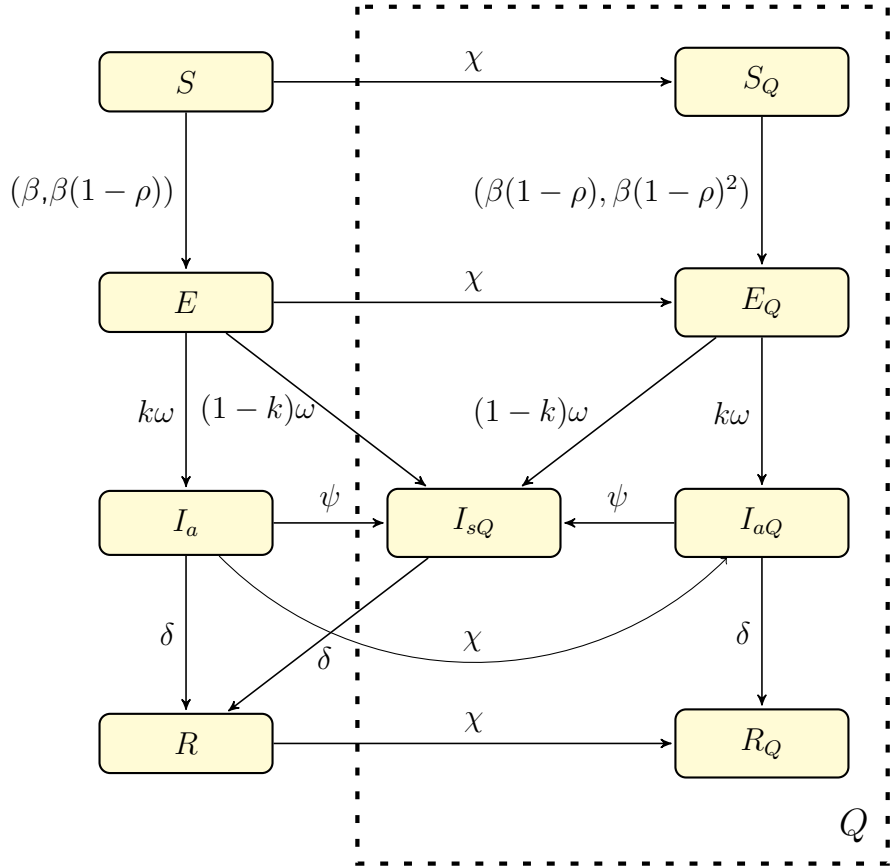

Several recent studies have reported that  
 the coronavirus can be transmitted during the incubation period before the first symptoms develop. Moreover, a significant portion of 
 individuals with COVID-19 lack symptoms \cite{asym}. 
These characteristics of  COVID-19 together with the state imposed quarantine  and testing motivate  the following extension of the SEIR model:
	\begin{equation}\label{model1}
	\begin{array}{l}
	\dot S= -\beta SI_a - \beta(1-\rho) S(I_{aQ}+ I_{sQ}) - \chi  I_{sQ} S,\\
	\dot S_Q= -\beta (1-\rho) S_QI_a -{ \beta}{(1-\rho)}^{2} S_Q(I_{aQ}+ I_{sQ}) +\chi I_{sQ} S,\\
	\dot E=  \beta SI_a  +\beta (1-\rho) S(I_{aQ}+I_{sQ})- \omega E-  \chi I_{sQ} E,\\
	\dot E_Q=  \beta (1-\rho) S_Q I_a +{ \beta}{(1-\rho)}^{2} S_Q(I_{aQ}+ I_{sQ}) - \omega E_Q+  \chi  I_{sQ}E   ,\\
	\dot I_a= k \omega E-  \chi I_{sQ}I_a -\delta I_a-\psi   I_a ,\\
	\dot I_{aQ}= k \omega E_Q+  \chi  I_{sQ}I_a   -\delta I_{aQ} -\psi   I_{aQ} ,\\
\dot I_{sQ}= (1-k) \omega (E+E_Q) + \psi (I_a+I_{aQ})  -\delta I_{sQ} ,\\
\dot R= \delta (I_a+ I_{sQ}) -\chi  I_{sQ} R,\\
\dot R_Q = \delta  I_{aQ}+\chi I_{sQ} R,
		\end{array}
	\end{equation}
%
see Figure \ref{Diagram1}. Here 
$S$  denotes the
density of susceptible individuals who are not quarantined, and $S_Q$ is the density of susceptible individuals at home 
quarantine.
A similar labeling convention is adopted for other 
groups of individuals, where 
the subscript $Q$ refers to the subpopulation at home quarantine.
It is assumed that the exposed individuals labeled $E$ and $E_Q$
have been infected but are not infectious yet, and
show negative test results if tested for the virus.
The compartment $I_{sQ}$ includes     
 infectious individuals with symptoms and
those infectious individuals without symptoms who tested 
positively for the virus (which  
corresponds to the published statistics such as the graphs provided by WHO  
\cite{who}). It is therefore the compartment of detected active cases.
The individuals from the compartment $I_{sQ}$ are all quarantined. 
Infectious asymptomatic individuals
labeled $I_a$ 
and $I_{aQ}$ show positive test results if tested for the virus,
in which case they are transferred to the compartment labeled $I_{sQ}$.
That is, testing an asymptomatic infectious non-quarantined individual results in their quarantining.
Upon recovery, the individuals from compartments $I_a$ and $I_{sQ}$ are recruited to the compartment $R$, and individuals from compartment $I_{aQ}$  are recruited to the compartment $R_Q$.
 We exclude reinfection, although a few cases have been reported. 
 The quarantined population is composed of the groups $S_Q,E_Q,I_{aQ},I_{sQ}, R_Q$.
 We assume a constant total population size, and we scale its density to unity, hence the 
equation for the density of the recovered 
quarantined population $R_Q=1-S-S_Q-E- E_Q-I_a-I_{aQ}-I_{sQ}-R$ 
is redundant. 

 We do not consider tests for antibodies which detect that a person had COVID-19 in the past but now is healthy.

Parameters $\beta$, $\omega$ and $\delta$  
have the same interpretation as for the SEIR model 
\eqref{SEIR}, i.e.~$\beta$ represents the  transmission rate, $\omega$ is the rate at which an exposed individual becomes infectious and $\delta$ is the recovery rate. 
The parameter $\psi$ denotes the 
testing 
rate, i.e. the  rate of 
detection and isolation of asymptomatic individuals. 
In addition, we assume that the total quarantining rate is $\chi  I_{sQ}  (S+I_a+E+R),$ where  the latter factor is the proportion of the non-quarantined population in the total population.
In particular, the total quarantining rate is proportional to $I_{sQ}$ and decreases with the increasing proportion of the total quarantined population.
Assuming indiscriminate quarantining, the rate of quarantining 
from a particular compartment is proportional to the density of the population associated with this compartment. For example, the rate of quarantining from the $S$-compartment is set to
$\chi  I_{sQ}  (S+I_a+E+R)  S/(S+I_a+E+R) = \chi  I_{sQ}  S$.
The model 
allows for the possibility of imperfect isolation with $\rho \in   [0, 1]$ representing the isolation effectiveness
as in \cite{11}. 
Finally, $k$ denotes  the 
fraction of individuas who do not develop symptoms 
when infected. Assuming $I_a(0)=I_{aQ}(0)=I_{sQ}(0)=0$, the ratio of the symptomatic population to the asymptomatic population remains $(1-k):k$ at all positive times.

 
The model exhibits 
a continuum of equilibrium states
for which  $I_a =  I_{sQ}=I_{aQ}= E_Q =E= 0$ and $S+S_Q + R+R_Q =1.$
 



\subsection{Model with additional features}
\label{Model3}
In this section, we consider a model with  a 
refined structure of  epidemiological compartments. 
The symptomatic and asymptomatic groups are considered separately 
because the transmission rate, the recovery rate and other epidemiological parameters can be different 
for these groups. We also 
include additional features of the disease such as the latent infectious period.

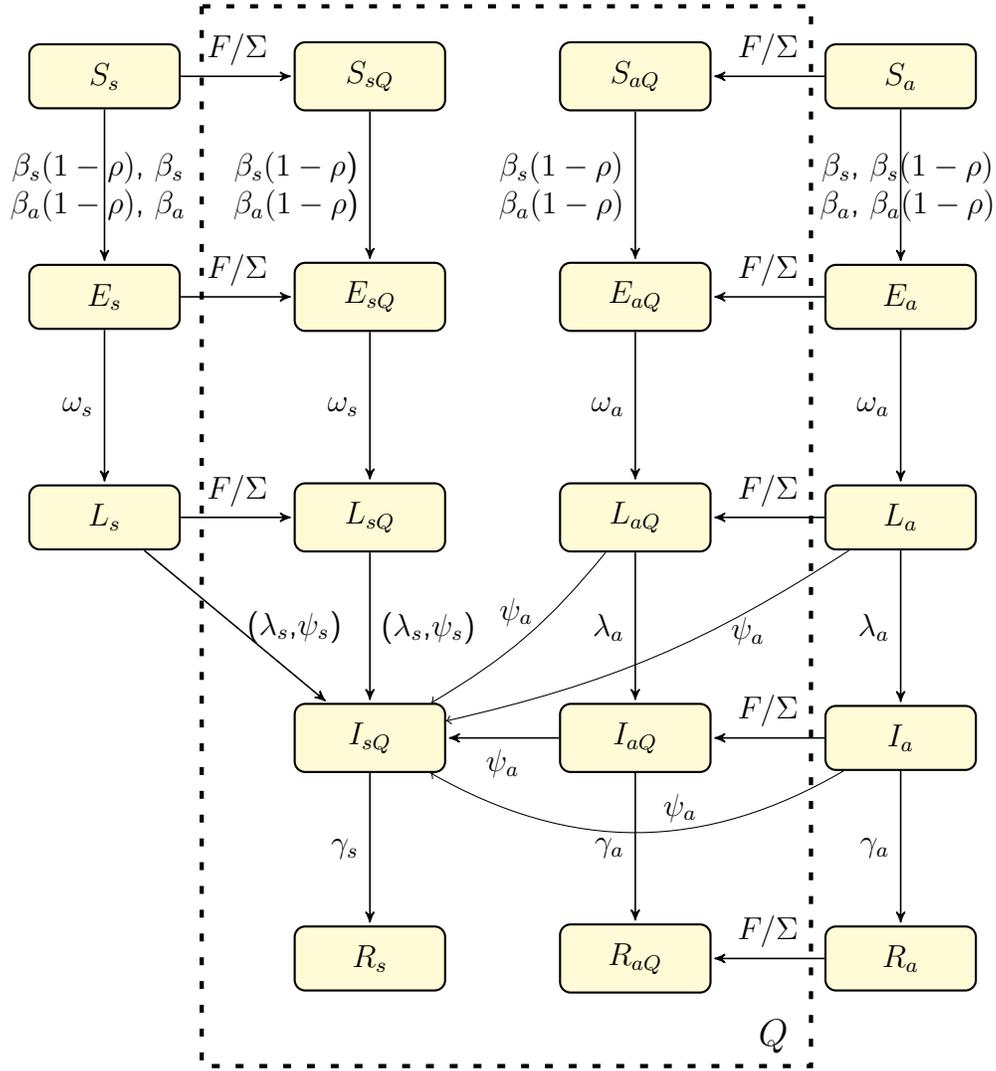
\begin{figure}
\begin{center}
\begin{tikzpicture}[
  font=\sffamily,
  every matrix/.style={ampersand replacement=\&,column sep=1.5cm,row sep=2cm},
  source/.style={draw,thick,rounded corners,fill=yellow!20,inner sep=.25cm,text width=1.5cm},
  to/.style={->,>=stealth',shorten >=1pt,semithick},
  every node/.style={align=center}]

  \matrix{
    \node[source](Ss){$S_s$};\&; \node[source](SsQ){$S_{sQ}$};  \&\   \node[source](SaQ){$S_{aQ}$};\&; \node[source](Sa){$S_{a}$}; \\
    \node[source](Es){$E_s$}; \&;\node[source](EsQ){$E_{sQ}$}; \&  \node[source](EaQ){$E_{aQ}$}; \&;\node[source](Ea){$E_{a}$};\\
      \node[source](Ls){$L_s$}; \&;\node[source](LsQ){$L_{sQ}$}; \&  \node[source](LaQ){$L_{aQ}$}; \&;\node[source](La){$L_{a}$};\\
     \& \node[source](IsQ){$I_{sQ}$}; \&; \node[source](IaQ){$I_{aQ}$}; \& \node[source](Ia){$I_{a}$};\\
     \& \node[source](Rs){$R_s$}; \&; \node[source](RaQ){$R_{aQ}$}; \& \node[source](Ra){$R_{a}$};\\
  };
  
    \draw[to] (Ss) -- node[midway,right] {{\hspace{-1.4cm}}$\beta_s(1-\rho)$,  $\beta_s$\\{\hspace{-1.4cm}}$\beta_a(1-\rho)$,  $\beta_a$}
    (Es);
        \draw[to] (SsQ) -- node[midway,left] {$\beta_s(1- \rho$)\\$\beta_a(1- \rho$)} (EsQ);
           \draw[to] (Sa) -- node[midway,right] {{\hspace{-1.22cm}}$\beta_s$, $\beta_s (1-\rho)$\\{\hspace{-1.22cm}}$\beta_a$, $\beta_a (1-\rho)$} (Ea);
        \draw[to] (SaQ) -- node[midway,left] {$\beta_s (1-\rho)$\\$\beta_a (1-\rho)$} (EaQ);
        
        \draw[to] (Es) -- node[midway,left] {$\omega_s$} (Ls);
        \draw[to] (EsQ) -- node[midway,left] {$\omega_s$} (LsQ);
           \draw[to] (Ea) -- node[midway,left] {$\omega_a$} (La);
        \draw[to] (EaQ) -- node[midway,left] {$\omega_a$} (LaQ);
        
     \draw[to] (Ss) -- node[midway,above] {$F/\Sigma$} (SsQ);
       \draw[to] (Sa) -- node[midway,above] {$F/\Sigma$} (SaQ);
       
          \draw[to] (Es) -- node[midway,above] {$F/\Sigma$} (EsQ);
       \draw[to] (Ea) -- node[midway,above] {$F/\Sigma$} (EaQ);
        \draw[to] (Ls) -- node[midway,above] {$F/\Sigma$} (LsQ);
       \draw[to] (La) -- node[midway,above] {$F/\Sigma$} (LaQ);
\draw[to] (Ia) -- node[midway,above] {$F/\Sigma$} (IaQ);
\draw[to] (La) -- node[midway,left] {$\lambda_a$} (Ia);
\draw[to] (Ls) -- node[midway,right] {($\lambda_s$,$\psi_s$)} (IsQ);
\draw[to] (LsQ) -- node[midway,right] {($\lambda_s$,$\psi_s$)} (IsQ);
               
     \draw[to] (LaQ) -- node[midway,left] {$\lambda_a$} (IaQ);
  \draw[to] (Ra) -- node[midway,above] {$F/\Sigma$} (RaQ);
  
  \draw[to] (IsQ) -- node[midway,left] {$\gamma_s$} (Rs);
    \draw[to] (Ia) -- node[midway,left] {$\gamma_a$} (Ra);
    \draw[to] (IaQ) -- node[midway,below] {$\psi_a$} (IsQ);
      \draw[to] (IaQ) -- node[midway,left] {$\gamma_a$} (RaQ);
      \path[->]
 (La) edge [bend left=10] node[midway,above] {{\hspace{2.3cm}} $\psi_a$} (IsQ);
 \path[->]
 (LaQ) edge [bend left=10] node[midway,above] {{\hspace{-0.4cm}} $\psi_a$} (IsQ);
\path[->]
 (Ia) edge [bend left=30] node[midway,above] {{\hspace{1cm}} $\psi_a$} (IsQ);
  \draw [loosely dashed,ultra thick] (-4,6.8) rectangle (4.1, -7.3);
  \node at (3.6,-6.9) {\large $Q$};
\end{tikzpicture}
\end{center}
 \caption{\label{Diagram2} Compartmental diagram of model \eqref{3}. Arrows show the flows between the  compartments
 with the corresponding rate parameters. The dashed rectangle includes the quarantined populations.  }
\end{figure}

We divide each compartment $S,S_Q,E,E_Q,R$ introduced in model \eqref{model1} into two compartments containing a subpopulation of individuals showing symptoms (labeled by the subscript $s$) and a subpopulation of asymptomatic individuals (labeled by the subscript $a$). 
For example, $S_{sQ}$ denotes the density of 
susceptible quarantined individuals who, if infected, will develop symptoms.
A similar indexing convention is adopted for the other groups.
In addition, we introduce the latent stage of infection
and the corresponding groups $L_s$, $L_{sQ}$, $L_a$, $L_{aQ}$ of symptomatic and asymptomatic individuals.
Figure \ref{Diagram2} presents the diagram of the compartments and the transitions between them for the following model:

		\begin{equation}\label{3}
	\begin{aligned}
	\dot S_s&= - S_s(\beta_s L_s+\beta_a(L_a+ I_a))  -  S_s (1-\rho)(\beta_s (L_{sQ}+ I_{sQ})\\&+\beta_a (L_{aQ}+I_{aQ})) -  \frac{S_sF}{\Sigma},\\
	\dot S_{sQ}&= -  S_{sQ}((1-\rho)\beta_s L_s+ (1-\rho)\beta_a (L_a+ I_a))  + \frac{S_sF}{\Sigma},\\
	\dot E_{s}&=S_s(\beta_s L_s+ \beta_a(L_a+ I_a))  + (1-\rho) S_s (\beta_s(L_{sQ}+ I_{sQ})\\&+\beta_a(L_{aQ}+I_{aQ}))   -\omega_s E_s - \frac{E_sF}{\Sigma}  ,\\
	\dot E_{sQ}&=  (1-\rho) S_{sQ}(\beta_s L_s+\beta_a(L_a+ I_a)) -\omega_s E_{sQ}+ \frac{E_sF}{\Sigma}, \\
	\dot L_{s}&=  \omega_s E_s   -(\lambda_s+ \psi_s) L_s - \frac{L_sF}{\Sigma}  ,\\
	\dot L_{sQ}&=  \omega_s  E_{sQ} -(\lambda_s+ \psi_s) L_{sQ} + \frac{L_sF}{\Sigma}, \\
	\dot I_{sQ}&=  (\lambda_s + \psi_s) L_s +   (\lambda_s + \psi_s) L_{sQ}-\gamma_s I_{sQ}  \\&+  \psi_a (L_a + L_{aQ}) +  \psi_a (I_a + I_{aQ}),\\
	\dot R_s &= \gamma_s I_{sQ},
	\end{aligned}
	\end{equation}
	\begin{equation*}
	\begin{aligned}
	\dot S_a&= - S_a(\beta_s L_s+\beta_a(L_a+ I_a))  -  (1-\rho) S_a(\beta_s(L_{sQ}+ I_{sQ})\\&+\beta_a(L_{aQ}+I_{aQ})) -  \frac{S_aF}{\Sigma},\\
	\dot S_{aQ}&= - (1-\rho) S_{aQ}(\beta_s L_s+\beta_a(L_a+ I_a))  + \frac{S_a F}{\Sigma},\\
	\dot E_{a}&=   S_a(\beta_s L_s+\beta_a(L_a+ I_a))  + (1-\rho) S_a (\beta_s (L_{sQ}+ I_{sQ})\\&+\beta_a(L_{aQ}+I_{aQ}))   -\omega_a E_a - \frac{E_a F}{\Sigma}  ,\\
	\dot E_{aQ}&=  (1-\rho) S_{aQ}(\beta_s L_s+\beta_a(L_a+ I_a)) -\omega_a E_{aQ} + \frac{E_a F}{\Sigma}, \\
	\dot L_{a}&=  \omega_a E_a  -\lambda_a L_a-\psi_a L_a - \frac{L_a F}{\Sigma}  ,\\
	\dot L_{aQ}&=   \omega_a E_{aQ}-\lambda_a L_{aQ}-\psi_a L_{aQ} + \frac{L_a F}{\Sigma}, \\
	\dot I_{a}&=  \lambda_a   L_{a} -  \psi_a I_a -\gamma_a I_{a}- \frac{I_aF}{\Sigma}  ,\\
\dot I_{aQ}&=  \lambda_a   L_{aQ}  -     \psi_a  I_{aQ}-\gamma_a I_{aQ} +\frac{I_a F}{\Sigma} ,\\
	\dot R_a &= \gamma_a I_{a}-\frac{R_a F}{\Sigma},\\
\dot R_{aQ} &= \gamma_a I_{aQ}+\frac{R_a F}{\Sigma},
\end{aligned}
	\end{equation*}
	where 
		\begin{equation}\label{3'}
\begin{aligned}	
		F &= \lambda_s L_{sQ} + \psi_s L_{sQ} + \psi_a(L_{aQ}+  I_{aQ}),\\
 \Sigma&= S_s  + S_a +E_s + E_a +L_s +L_a+ I_a+ R_a.
\end{aligned}
\end{equation}
 The parameters $\beta_s$, $\beta_a$, $\omega_s$, $\omega_a$, $\psi_s$, $\psi_a$ and $\rho$ have the same meaning as the corresponding parameters in model \eqref{model1}.
 The recovery rates $\gamma_s$ and $\gamma_a$ are higher than their counterpart $\delta$ in model \eqref{model1} because the latent phase
 of the average duration $1/\lambda_s$ and $1/\lambda_a$  for the symptomatic and asymptomatic groups, respectively, is introduced in model \eqref{3}.

Let us note a few differences between models \eqref{model1} and \eqref{3}.
In the latter model,  the initial stage of infection (represented by the groups labeled $E$)
when the individuals are not yet infectious is followed by 
the latent infectious incubation stage (represented by the groups labeled $L$)\footnote{The latency stage can be similarly introduced into model \eqref{model1}.}. 
This allows us to account for the fact that individuals can become contagious and spread the disease before developing symptoms, which is an important characteristic of COVID-19 although the latency period is relatively short \cite{Incubation}.  
Further, dividing the compartments as in model \eqref{3}, one can use different parameters for symptomatic and asymptomatic groups (i.e., groups labeled $s$ and $a$), and thus emphasize that individuals without symptoms (or with mild symptoms)
and the symptomatic cases can possibly show different epidemiological characteristics. For instance, 
a difference in the transmission rate and the average duration of the infectious and latent periods for asymptomatic and symptomatic individuals
can be accounted for if needs be.
Moreover,  one can consider different values of the testing rate for different compartments reflecting the fact that the disease is harder to detect in asymptomatic individuals. 
 
Models \eqref{model1} and \eqref{3} use somewhat different assumptions about the recruitment to quarantine.  Model \eqref{model1} postulates a constant rate
of testing with the associated isolation of positively testing individuals and, additionally, a constant rate of isolation unrelated to testing.
In model \eqref{3}, the proportion of the total quarantined population,  $\Sigma_Q= S_{sQ} + E_{sQ} + L_ {sQ} + S_{aQ} +E_{aQ}+ L_ {aQ} +I_{aQ}+R_{aQ},$
is conserved as we assume that $\Sigma_Q$ is maintained at the same level enforced by the government interventions. 
The relation $\dot \Sigma_Q=0$ is warranted by a feedback loop, which relates the recruitment rate $F$ to the instantaneous values of the population densities by formula \eqref{3'}.



 Contrary to  model \eqref{model1}, and as a simplification, we assume zero transmission
from the quarantined individuals because it is relatively small
 (however, the corresponding terms proportional to $(1-\rho)^2$
 can also be included in \eqref{3} as in model \eqref{model1}).

Model \eqref{3} also exhibits  
a continuum of equlibrium states 
for which  $I_a =L_a=E_a=  R_s=R_a= I_{sQ}=I_{aQ}=E_a= E_{aQ} =E_s=E_{sQ}= S_s=L_s=L_{sQ} = S_a= 0$ and $S_{aQ} +S_{sQ}+ R_{aQ} =1.$
 
\subsection{
Epidemiological parameters}
\label{parameters}

The values  of the epidemiological parameters used in the numerical simulations presented below are based on the average characteristics of COVID-19 published so far.     

One key parameter 
is the average duration of the infectious period. The period of infectiousness (or period of communicability) is defined as \emph{the time interval during which an infectious agent may be transferred directly or indirectly from an
infected person to another person} \cite{cdc}. According to  \cite{Incubation}, COVID-19 can be transmitted by an infected individual before the  symptoms develop. 
While in model \eqref{model1} the average infectious period is the reciprocal of the parameter $\delta$, in model \eqref{3} the infectious period is divided into two stages with the associated populations $L $ and $I $ and is represented by the parameters $\lambda$ and $\gamma.$ 
In order to keep our results for models \eqref{model1}  and \eqref{3}  comparable,
we assume that $1/\delta=1/\lambda + 1/\gamma$.

The period of infectiousness is difficult to estimate directly. 
One of a few current studies including 
nine patients with mild symptoms 
showed that all patients were non-infectious
on the eighth day after the first symptoms developed \cite{nature}. 

Another possibility to estimate the parameter $\delta$ is through the estimation of the so-called \emph{serial interval}. The serial interval of an infectious disease represents the duration between symptom onset of a primary case and symptom onset of its secondary cases.
According to \cite{Incubation}, the mean incubation period of COVID-19 is 6 days. In \cite{SI2}, the serial interval of COVID-19 was estimated as the weighted mean of the published parameters and 
described by the gamma distribution with the mean serial interval of 4.56 days (credible interval $(2.54,7.36)$) and standard deviation 4.53 days (credible interval $(4.17-5.05)$).

With a serial interval shorter than the average incubation period ($4.56<6$), we expect that a significant number of transmissions occur before the index case has symptoms. Following \cite{kucharski}, the  period of infectiousness is assumed to be 5.5 days, 2 days before  and 3.5 days after symptom onset.  
This data 
translates into the values $\lambda=1/2,$ $\gamma=1/3.5,$ $\delta=1/5.5,$  $\omega =1/4 $ 
of model parameters. 
The study \cite{nature} indicates that  the period of infectiousness may last longer for patients with symptoms (7 days),  hence we 
 use different values of $\gamma$ for the compartments labeled $s$ and $a$ in 
 Section \ref{sec5.3}.

Another key parameter
is the transmission rate $\beta$, which can be expressed as the product of the average number of daily contacts which a susceptible individual has with infected individuals and the probability of transmission during each contact. The value of this  parameter, which is not directly observable, is inferred from the estimation of the basic reproduction number $R_0$. 
The
value of $R_0$ for COVID-19 has been estimated within the range of $2.24\le R_0\le  3.58$ \cite{W}. Taking into  account the theoretical relation between $R_0$ and  $\beta$ in SEIR model,  $\beta = \delta R_0$, and the setup of the other model parameters, 
in our numerical analysis we assume $\beta=0.5.$ In model \eqref{3}, while we consider different transmission rates for symptomatic and asymptomatic individuals, we explore the range of $\beta$ from the interval $[0.1, 0,6].$

The coefficient $(1-\rho)\in [0,1]$  
measures the decrease in the transmission rate due to isolation.
For illustrative purposes, we set $\rho=0.5$. 
%
%
The value of the parameter $k$, which measures the proportion of asymptomatic individuals among the infected population, is set to $0.5$ in accordance with the estimate from \cite{who2}.









 
\section{
Results} \label{sec3}

The current epidemiological situation is raising many questions about the dynamics of the virus spread and the effectiveness of the imposed quarantine policies. The  proposed  models \eqref{model1} and \eqref{3} aim to  explore some plausible scenarios depending on the size of initial quarantine, gradual quarantining rate and the testing rate.

If not explicitly stated otherwise, the values of epidemiological parameters used in the simulations presented below are as follows: $\beta=0.5$, $\rho=0.5$, $\omega=\omega_s=\omega_a=1/4$, $\delta=1/5.5$, $\lambda_s=\lambda_a=1/2,$  $\gamma_s=\gamma_a=1/3.5$ and $k=0.5.$ 
The initial conditions are chosen in accordance with the real situation in the city where the population size $N$ is approximately $500, 000$. 
If not stated otherwise,  simulations of model \eqref{model1} are initiated with $E(0)=2/N$  and simulation of model \eqref{3}  with $E_s(0)=E_a(0)=1/N$. Moreover, we assume that originally the rest of the population is susceptible,
i.e. the initial size of all the groups except $E$ and $S$ is zero.
The values of parameters $\chi$, $\psi$, $\psi_s$ and $\psi_a$ will be specified for each particular simulation. 


\begin{figure}[h]
\begin{subfigure}[b]{0.5\textwidth}
\begin{center}
\includegraphics[width=7cm]{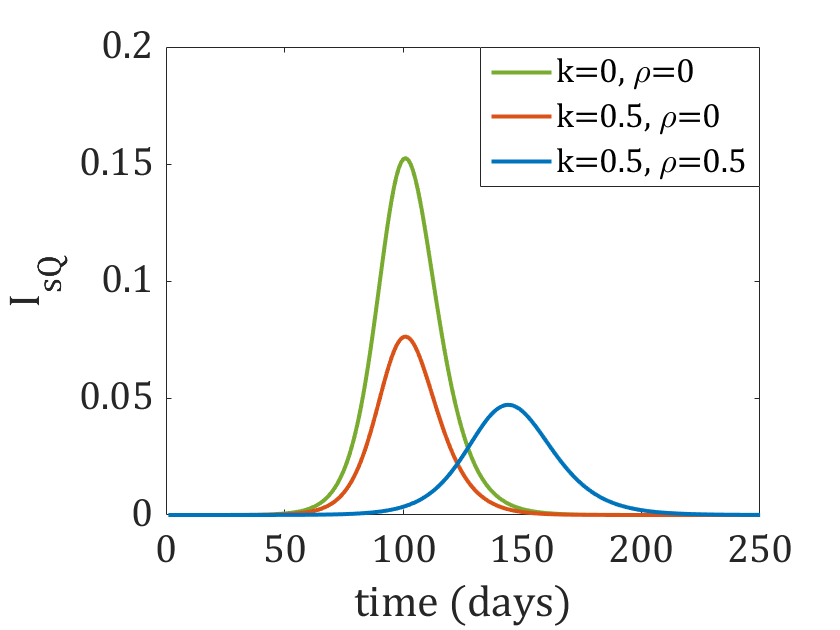}
\caption{}
\end{center}
\end{subfigure}
~
\begin{subfigure}[b]{0.5\textwidth}
\begin{center}
\includegraphics[width=7cm]{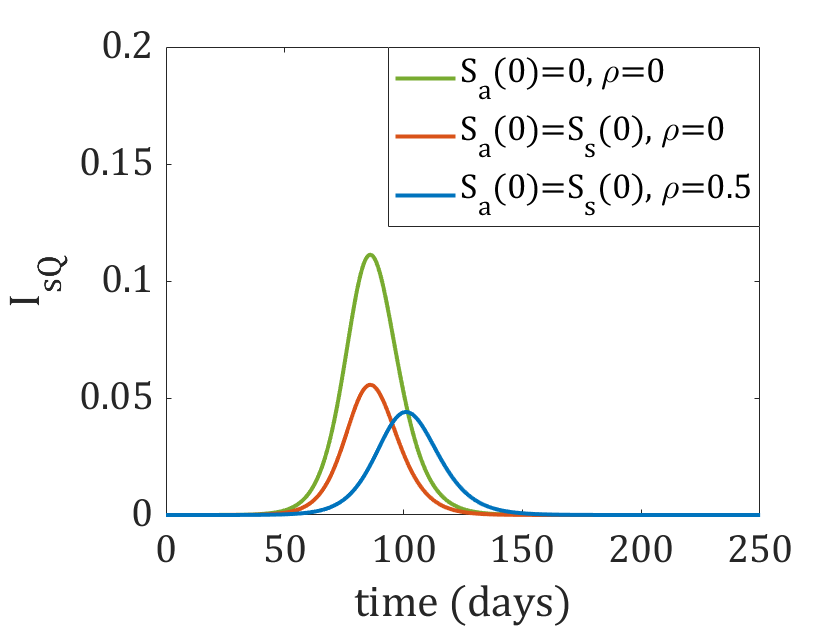}
\caption{}
\end{center}
\end{subfigure}
\begin{subfigure}[b]{0.5\textwidth}
\begin{center}
\includegraphics[width=7cm]{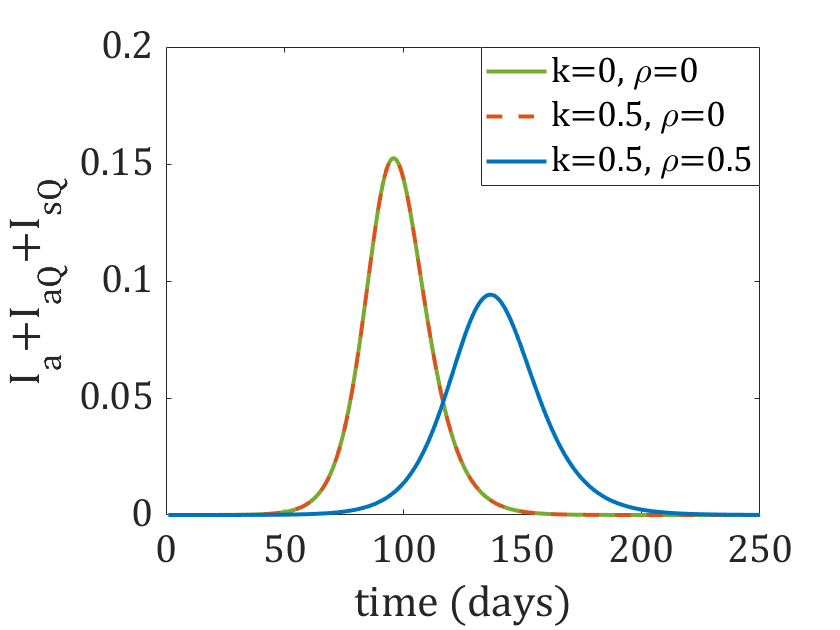}
\caption{}
\end{center}
\end{subfigure}
~
\begin{subfigure}[b]{0.5\textwidth}
\begin{center}
\includegraphics[width=7cm]{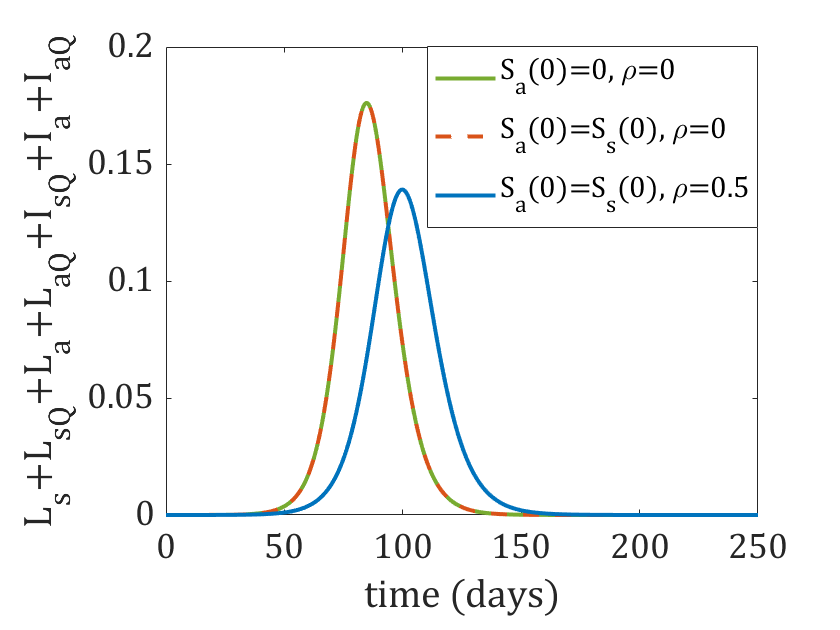} 
\caption{}
\end{center}
\end{subfigure}
\caption{\label{por} (a, c) Comparison of   the graphs of $I_{sQ}$ and $I_a+I_{aQ}+I_{sQ}$ for  model \eqref{SEIR} and  model  \eqref{model1} with parameters $\chi=0$, $\psi=0$ and initial conditions  $E(0)=2/N$.
(b, d) $I_{sQ}$ and the total infected population ($L_s+L_{sQ}+
L_a+L_{aQ}+I_a+I_{aQ}+I_{sQ}$)  for model \eqref{3} with $\psi_s=\psi_a=0$ and initial conditions $E_s(0)=2/N$  and $E_a(0)=0$ (green) and 
$E_s(0)=E_a(0)=1/N$ (red and blue). 
}
\end{figure}

We  first compare the classical SEIR model \eqref{SEIR} and model \eqref{model1}, see Figure \ref{por} (a). 
Model \eqref{model1} reduces to
\eqref{SEIR} in the special case when $k=0$, $\rho = 0$, $\chi= 0$, $\psi= 0$.  The green curve in Figure \ref{por}    (a) depicts the time evolution of the infected population for this case. Let us note that the  $I_s$ compartment in model \eqref{model1} fully corresponds to the  $I$ compartment in  \eqref{SEIR} under this parameters setup.
Taking $k=0.5$ with other parameters unchanged in model \eqref{model1},
we observe the decrease of the infected population by a half because half of the infected individuals don't show symptoms\,---\,instead of being included in the $I_s$ they are now logged in the compartment $I_a$. 
Setting $\rho = 0.5$ leads to a slower infection spread due to quarantine. 
Similar results for model \eqref{3} are shown in Figure \ref{por} (b), where instead of changing the parameter $k$ we 
change the initial proportion of symptomatic and asymptomatic groups. 

Let us note that typical COVID-19 statistics such as graphs provided by WHO \cite{who} record confirmed cases only. However, according to some sources including \cite{who2}, the data suggest that up to 80\% (we assume 50\% in simulations) of infections can be mild or asymptomatic. Therefore, one can expect that the  propagation of the disease by unregistered cases might play a significant role in its spread, possibly even more significant than the symptomatic cases because the latter are usually quarantined and therefore spread the virus less. 
It remains a controversial topic whether symptomatic and asymptomatic cases are equally infectious and in particular 
whether symptomatic and asymptomatic individuals produce the same numbers of antibodies, which would suggest that they spread the disease approximately equally. 
We assume equal transmission rates from symptomatic and asymptomatic individuals in model \eqref{model1} and allow for different transmission rates in model \eqref{3}. 

Figure \ref{por} shows that the proportion of asymptomatic individuals in the population can have a significant impact on the height of the maximum of $I_{sQ}$ (green versus  red and blue curves),   i.e.\ the number of infected symptomatic individuals at the peak of the epidemic. As a result, countries  with different proportions of symptomatic versus asymptomatic individuals (e.g.\ due to different age structure of the population or for other immunity reasons) might show different dynamics of the epidemic, although other conditions are similar. 
Introducing a reduction of the transmission due to quarantine by setting $\rho= 0.5$ lowers the infection peak and simultaneously prolongs the epidemic (red versus blue curve). This prolongation effect is more pronounced in model \eqref{model1}.
In model \eqref{3},
which extends models \eqref{SEIR} and \eqref{model1} by introducing the $L$ compartments, we observe a lower peak
of reported cases,
$I_{sQ}$, because individuals at the latent phase of the disease are not counted to this compartment.
On the other hand, the total number of infected individuals is higher than in model \eqref{model1}
because the infected individuals are not quarantined during the latent phase.
Also, the epidemic is shorter compared to model  \eqref{model1}.





\subsection{Response of the models to variations of the quarantine size}
\label{quarantine}
The most effective way to slow down the spread of a disease which is not vaccine preventable, is probably to impose isolation and  quarantine  on the population or a selected subpopulation. We assume that symptomatic and positively testing individuals are isolated (either at a hospital when symptoms are severe, or at home). Further, the models assume home quarantining of part of the untested population, including healthy and asymptomatic individuals who stay at home due to various state/business restrictions (such as school closures etc.).  We assume that 
these interventions result in a decreased transmission rate.  

In both   models \eqref{model1} and \eqref{3} home quarantine can be controlled in two possible ways. 
One is through the initial conditions. The size of the quarantine  is adjusted by setting $S_Q(0) > 0$ and $S_{sQ}(0)+S_{aQ}(0)>0$ for models \eqref{model1} and \eqref{3}, respectively. This type of quarantine is further referred to as \emph{abrupt} quarantining. 

\begin{figure}[h!]
\begin{subfigure}[h]{0.5\textwidth}
 \begin{center}
\includegraphics[width=7cm]{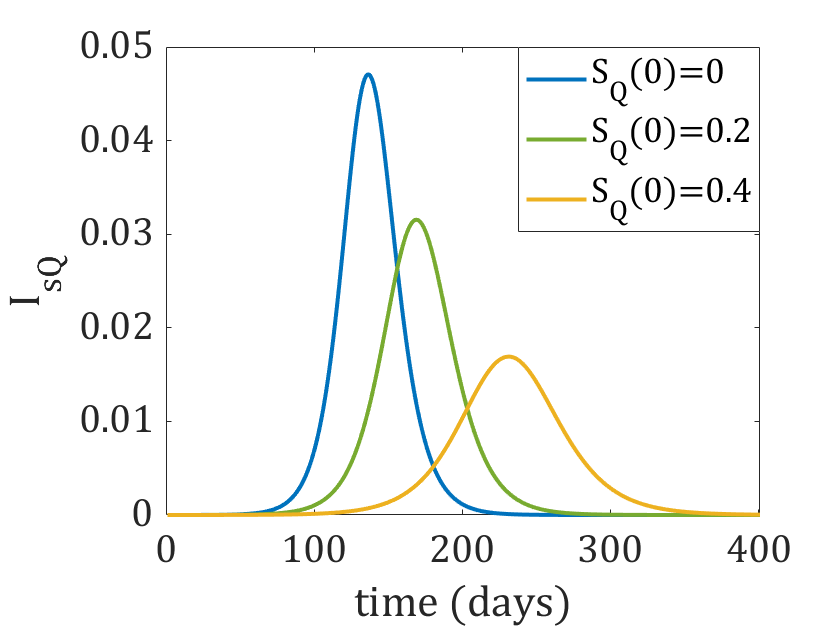}
\caption{}
\end{center}
\end{subfigure}%
~
\begin{subfigure}[h]{0.5\textwidth}
\begin{center}
\includegraphics[width=7cm]{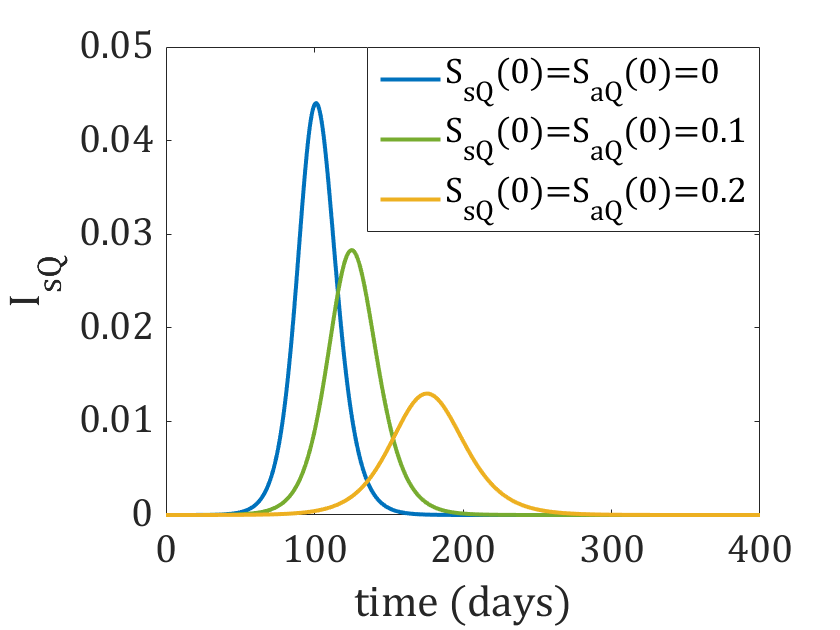}
\caption{}
\end{center}
\end{subfigure}

\begin{subfigure}[h]{0.5\textwidth}
 \begin{center}
\includegraphics[width=7cm]{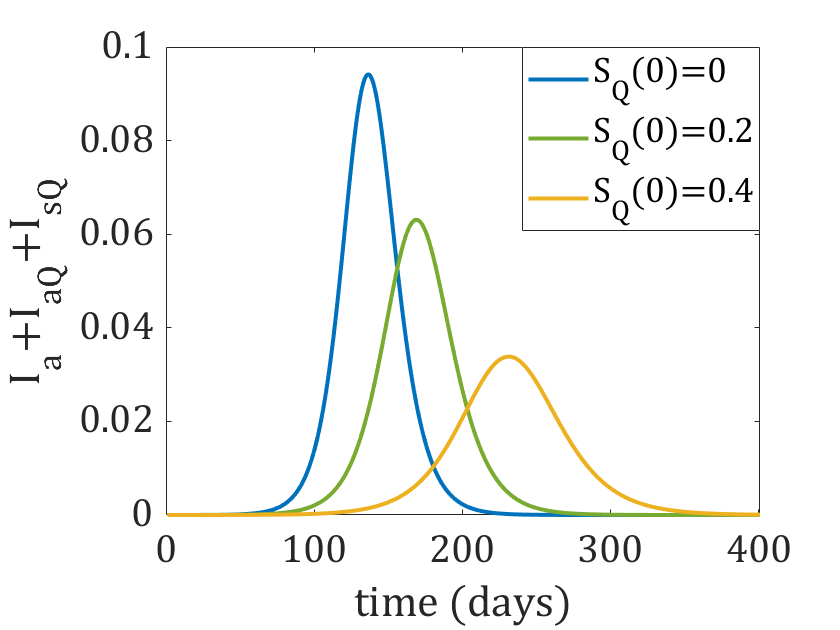}
\caption{}
\end{center}
\end{subfigure}%
~
\begin{subfigure}[h]{0.5\textwidth}
\begin{center}
\includegraphics[width=7cm]{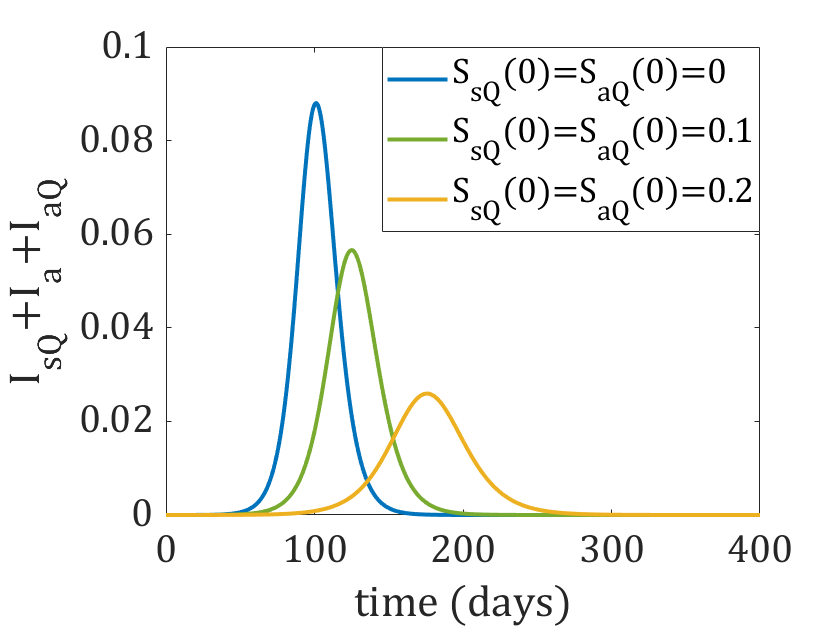}
\caption{}
\end{center}
\end{subfigure}

\caption{Effect of the quarantine size on $I_{sQ}$ and $I_a+I_{aQ}+I_{sQ}$. (a,c)  Model \eqref{model1} with $k=0.5$, $\rho=0.5$, $\psi=0$, $\chi=0$ and $E(0)=2/N$.
(b,d) Model \eqref{3}  with $\psi_s=\psi_a=0$ and $E_s(0)=E_a(0)=1/N$, $S_s(0)=S_a(0)$. 
\label{SQ}}
\end{figure}

Figure \ref{SQ} shows how the the number of infected symptomatic individuals depends on the  initial quarantine size.  As expected, the infection peak of $I_{sQ}$ lowers with increasing $S_Q(0)$ ($S_{sQ}(0)+S_{aQ}(0)$, respectively), but simultaneously the duration of the epidemic increases. 
This observation is similar for both models. 
One can see that 
the duration of the epidemic is shorter and the height of the infection peak is lower for model \eqref{3},
which can be attributed to the presence of the latent phase $L$.



Another way to control the home quarantine size in model \eqref{model1} is 
by using a positive
isolation rate $\chi$. Below this strategy is referred as the \emph{gradual} quarantining. 
The setup allows $\chi$ to depend on time or to be controlled by the phase variables such as $I$ through a feedback loop but, for simplicity, we assume $\chi$ to be constant.

\begin{figure}[h]
\begin{subfigure}[h]{0.5\textwidth}
\begin{center}
\includegraphics[width=7cm]{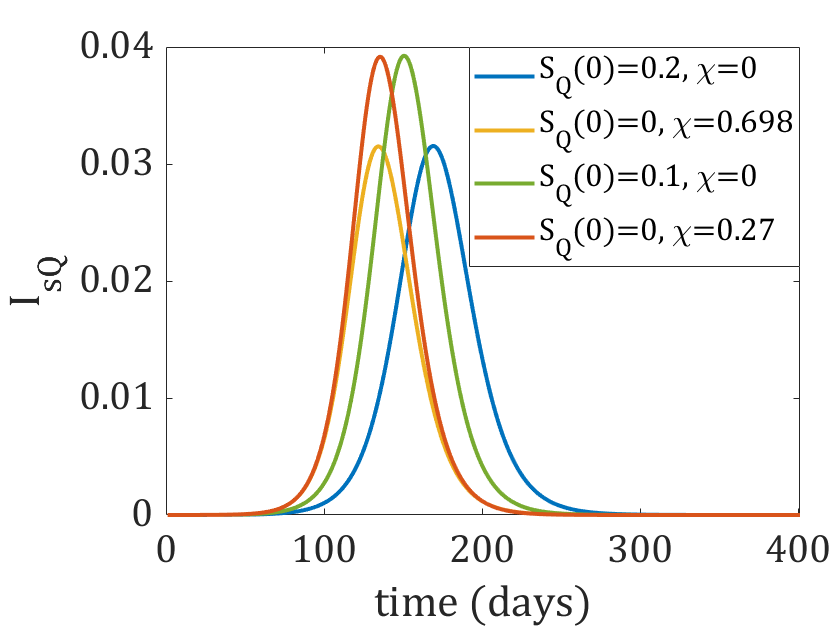}
\caption{}
\end{center}
\end{subfigure}%
~
\begin{subfigure}[h]{0.5\textwidth}
\begin{center}
\includegraphics[width=7cm]{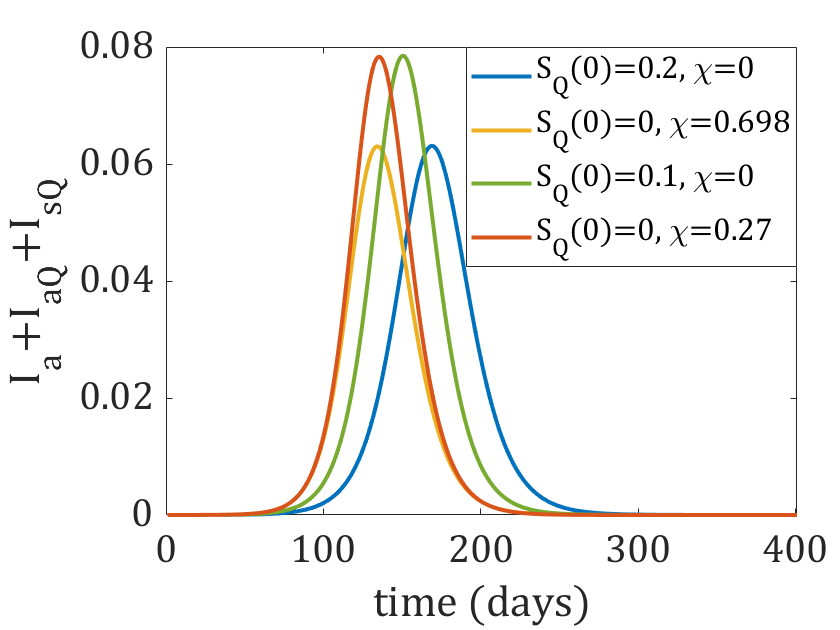}
\caption{}
\end{center}
\end{subfigure}
\begin{subfigure}[h]{0.5\textwidth}
\begin{center}
\includegraphics[width=7cm]{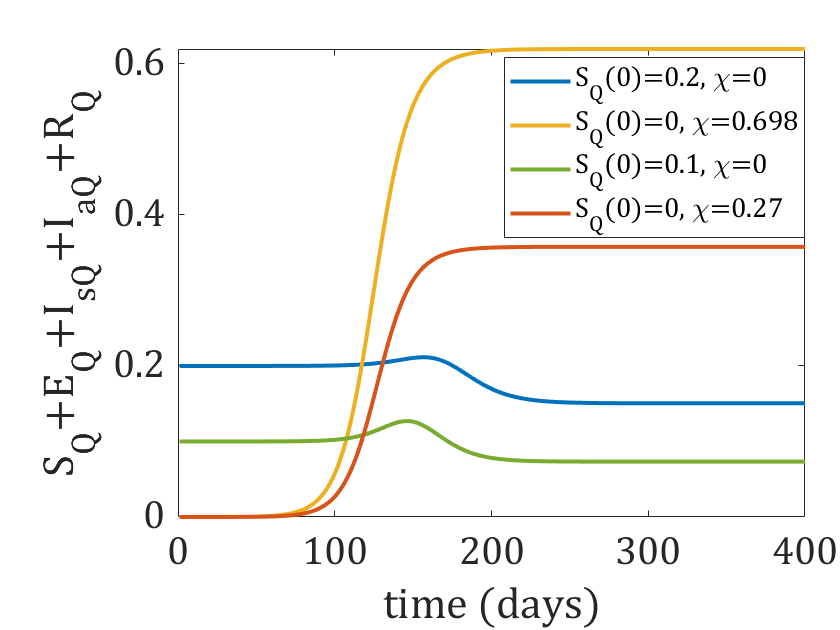}
\caption{}
\end{center}
\end{subfigure}%

\caption{Comparison of the effect of abrupt and gradual quarantining on dynamics of (a) $I_{sQ}$ population; (b) $I_{a}+I_{aQ}+I_{sQ}$ population; and (c) quarantined population in model \eqref{model1}. 
\label{two_quarantine}}
\end{figure}

Figure \ref{two_quarantine} compares the above two approaches to controlling the quarantine size in model \eqref{model1}. The same infection peak and a similar profile of the function $I_{sQ}(t)$ can be achieved either with nonzero $\chi$ or with nonzero $S_Q(0),$ see Figure \ref{two_quarantine} (a).  The total quarantined population 
is twice larger than $I_{sQ}$ because equal proportions of symptomatic and asymptomatic individuals are assumed, see Figure \ref{two_quarantine} (b). 
Figure \ref{two_quarantine} (c) shows that the total quarantined population is initially smaller for the gradual quarantining
strategy than for the abrupt quarantine strategy but eventually has to exceed the latter to ensure the same infection peak height.
The peak appears  earlier with gradual quarantining, and 
this effect becomes more pronounced with increasing  quarantine size (blue vs. yellow curve).






\begin{figure}[h]
\begin{subfigure}[h]{0.5\textwidth}
\begin{center}
\includegraphics[width=7cm]{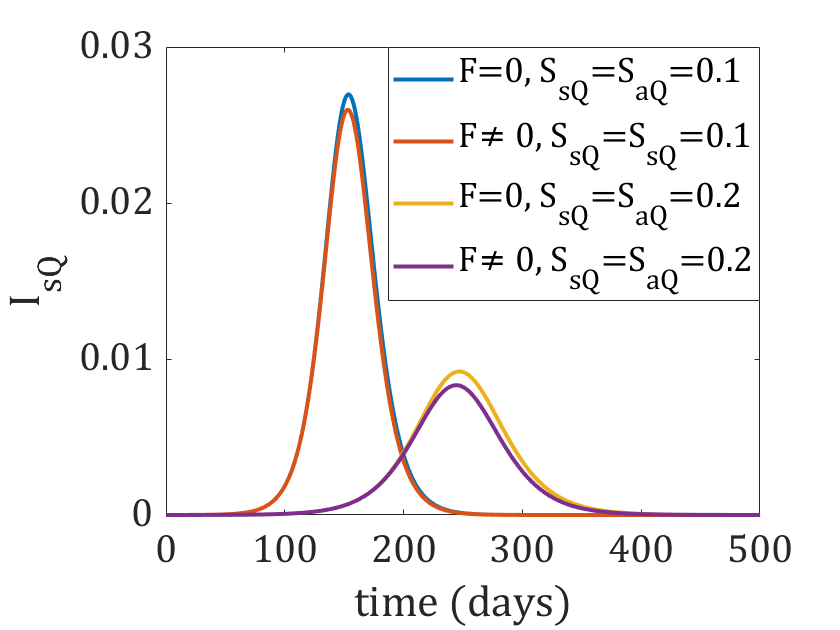}
\caption{}
\end{center}
\end{subfigure}%
~
\begin{subfigure}[h]{0.5\textwidth}
\begin{center}
\includegraphics[width=7cm]{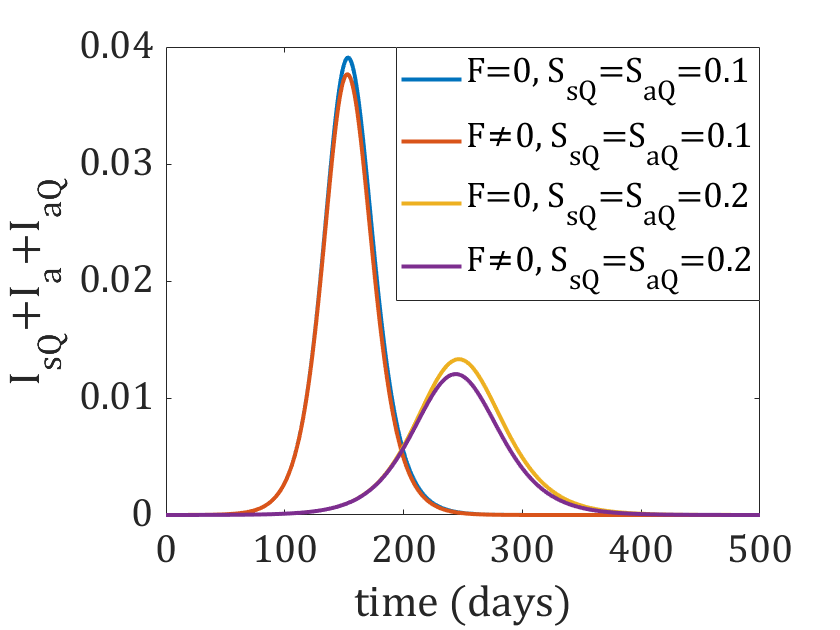}
\caption{}
\end{center}
\end{subfigure}
\caption{Effect of the conservation condition on (a) $I_{sQ}$  and (b) $I_{sQ}+I_{a}+I_{aQ}$ in model \eqref{3} for two different sizes of the initial quarantine and    $\psi=0.1$. The value of $F$ is defined by \eqref{3'} (the size of the quarantine is conserved) for the red and purple plots and set to zero  (conservation is violated) for the blue and yellow plots. \label{Fvalue}}
\end{figure}

In model \eqref{3},
 the quantity
$\Sigma= S_{sQ} + L_{sQ} + S_{aQ} +L_{aQ}+ I_{aQ}+R_{aQ}$, which represents the total quarantined population, is conserved
(more realistically, the quarantined population cannot exceed a certain value due to economic reasons, the necessity of uninterrupted operations of essential services etc.). 
Let us note that the isolation rate $F$ is actually a small quantity depending on the testing rate and the number of people in quarantine.   
The effect of the conservation condition represented by the value  $F$  is illustrated in Figure \ref{Fvalue}. 
Here we assume that initially $S_{sQ} (0) + S_{aQ}(0) = 0.2$ for the blue and red curves and $S_{sQ} (0) + S_{aQ}(0) = 0.4$ for the yellow and purple curves and $\psi = 0.1$. 
The blue curve represents the plot of $I_{sQ}$ for model \eqref{3} with $F,\Sigma$ defined by \eqref{3'},
while the red curve is obtained with  $F=0$ for the same parameters.
The purple and yellow curves were obtained similarly
with a different initial quarantine size.  
For all the selected setups of parameters, $F$ does not seem to have a significant effect on dynamics of the number of reported cases 
and the total infected population.


%

\subsection{Response of the models to variations of the testing rate}

In both models \eqref{model1} and \eqref{3}, by testing we understand the detection and isolation of infectious asymptomatic individuals. In model \eqref{3}, this category includes also   individuals during the latent period. We assume that symptomatic cases are confirmed as positive and  isolated (quarantined), hence they all belong to the $I_{sQ}$ compartment.
The parameter  $\psi$  can be interpreted as the rate of successful detection of infectious asymptomatic individuals.  We assume that this rate is the same for all  compartments in which testing can detect infectious individuals. Newly detected cases are recruited from  the compartments  $I_a$ and $I_{aQ}$ in model \eqref{model1} and from the compartments $L_a$, $L_{aQ}$, $L_s$, $L_{sQ}$, $I_a$,  $I_{aQ}$ in model \eqref{3}, respectively, to  $I_{sQ}$. 
Therefore the compartment $I_{sQ}$ is interpreted  as the compartment of \emph{ confirmed active} cases. 

\begin{figure}[h!]
\begin{subfigure}[b]{0.5\textwidth}
\begin{center}
\includegraphics[width=7cm]{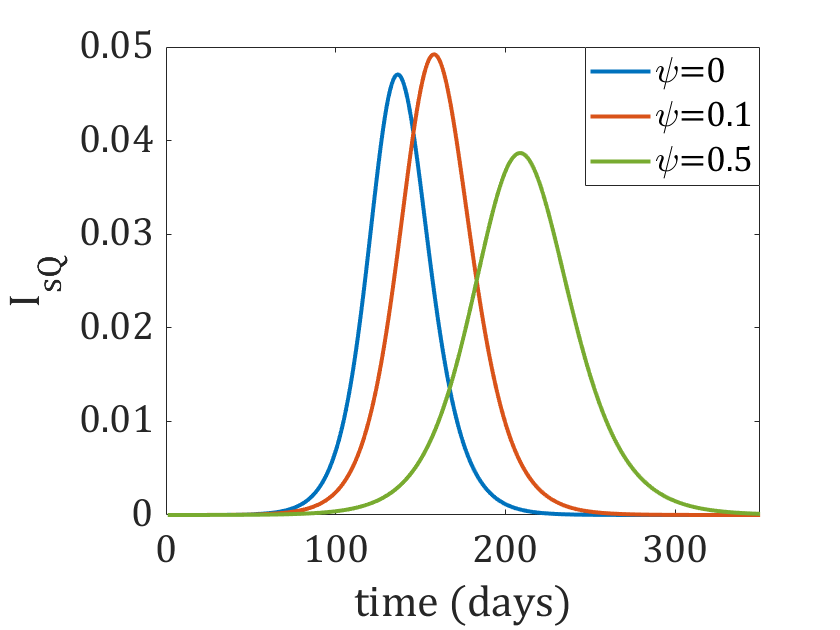}
\caption{}
\end{center}
\end{subfigure}%
~
\begin{subfigure}[b]{0.5\textwidth}
\begin{center}
\includegraphics[width=7cm]{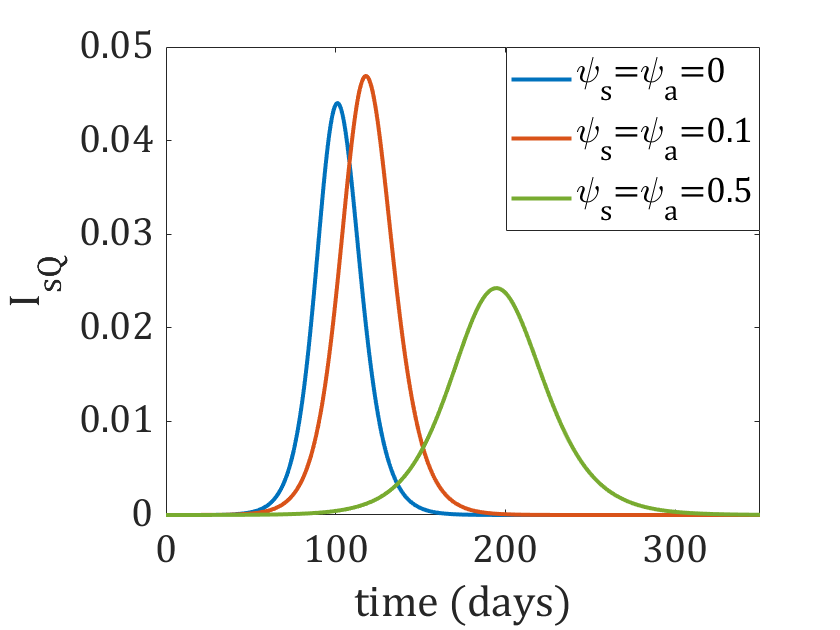}
\caption{}
\end{center}
\end{subfigure}



\begin{subfigure}[b]{0.5\textwidth}
\begin{center}
\includegraphics[width=7cm]{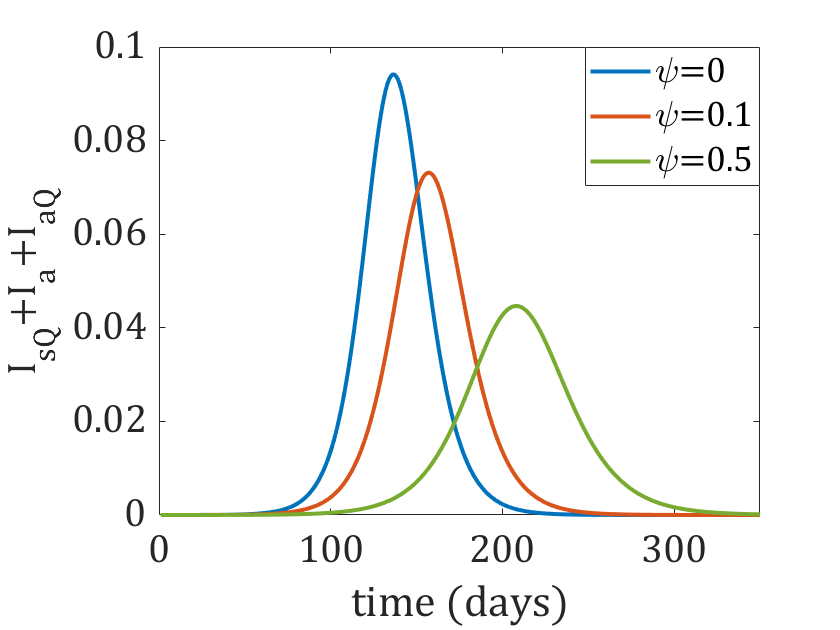}
\caption{}
\end{center}
\end{subfigure}%
~
\begin{subfigure}[b]{0.5\textwidth}
\begin{center}
\includegraphics[width=7cm]{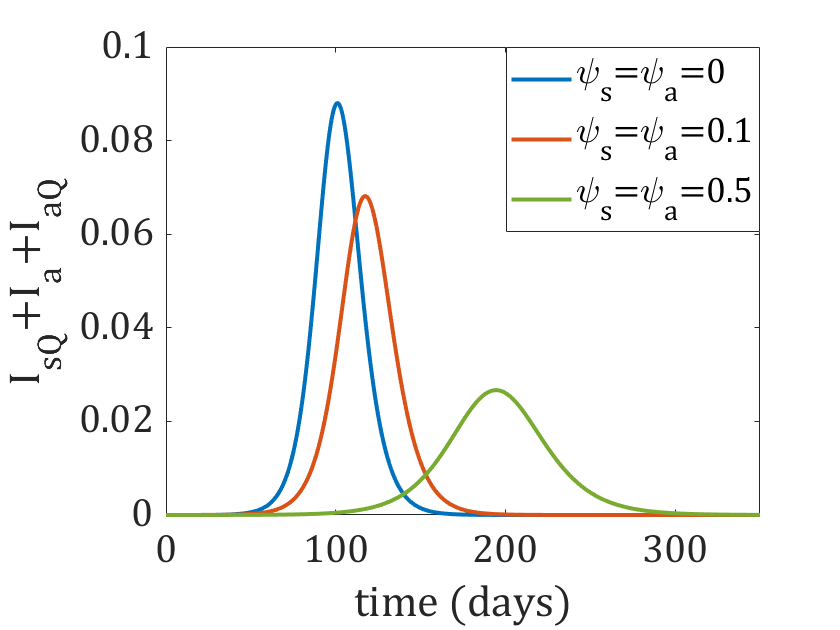}
\caption{}
\end{center}
\end{subfigure}
\caption{ Effect of testing on dynamics of $I_{sQ}$  and $I_a+I_{aQ}+I_{sQ}$. (a, c) Model \eqref{model1} with $\chi=0$ and $E_Q(0)=2/N$. 
(b, d) Model \eqref{3} with  $E_{s}(0)= E_a(0)=1/N$. 
\label{PSI1} }

\end{figure}

As our simulations show,  
 the time interval from the advent of the epidemic to the infection peak and the total duration of the epidemic increase with the increasing testing rate $\psi$, see Figure~\ref{PSI1} (a) for model \eqref{model1}  and (b) for model \eqref{3}, respectively. Further,   a significant increase in the testing rate decreases the epidemic peak.
This effect is more pronounced for model \eqref{3}.
However, for a relatively low testing rate, for instance if  $\psi = 0.1$, a small increase in the number of confirmed active cases is observed, compared to the case when $\psi=0$ (the red curve versus the blue curve in  Figure~\ref{PSI1}). Under low detection rate,  the  decrease of the spread of infection due to quarantining of positively tested asymptomatic individuals can be masked by
the apparent increase  of the $I_{sQ}$ population due to newly detected cases. However, the total number of cases (i.e. $I_{sQ} + I_{aQ} + I_a$ for model \eqref{model1} and $I_{sQ}+L_s+L_{sQ}+L_a+L_{aQ}+ I_a + I_{aQ}$  for model \eqref{3}, respectively) always decreases with increasing testing rate, as Figure \ref{PSI1} (c,d) shows.

\subsection{The effect of the difference of parameters of symptomatic and asymptomatic groups 
} \label{sec5.3}



The composition of model \eqref{3} allows us to use different parameters for symptomatic and asymptomatic cases and thus account for 
the possibility that the COVID-19 illness can take different course for these two groups. 
In principle, symptomatic and asymptomatic groups can differ in any  model parameters. 
 We tested the sensitivity of the $I_{sQ}$ and $I_a + I_{aQ} + I_{sQ}$ peak heights to the average length of the infectious stage $1/\gamma_s$ and $1/\gamma_a$ in the symptomatic and asymptomatic groups, respectively; the transmission rates $\beta_s$ and $\beta_a$; and the testing based quarantining rates $\psi_s$ and $\psi_a$. 
 
\begin{figure}[h!]
\begin{subfigure}[b]{0.5\textwidth}
\begin{center}
\includegraphics[width=7cm]{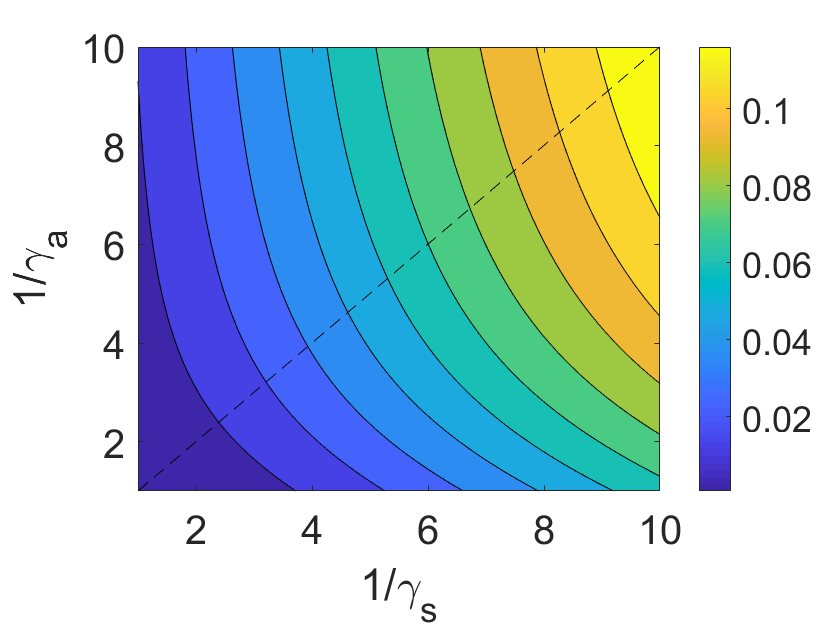}
\caption{}
\end{center}
\end{subfigure}%
~
\begin{subfigure}[b]{0.5\textwidth}
\begin{center}
\includegraphics[width=7cm]{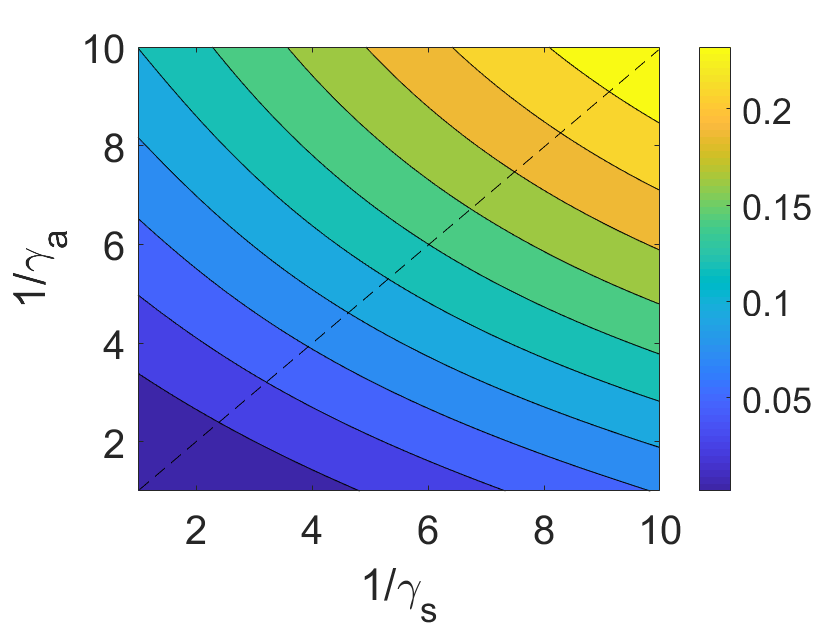}
\caption{}
\end{center}
\end{subfigure}
\begin{subfigure}[b]{0.5\textwidth}
\begin{center}
\includegraphics[width=7cm]{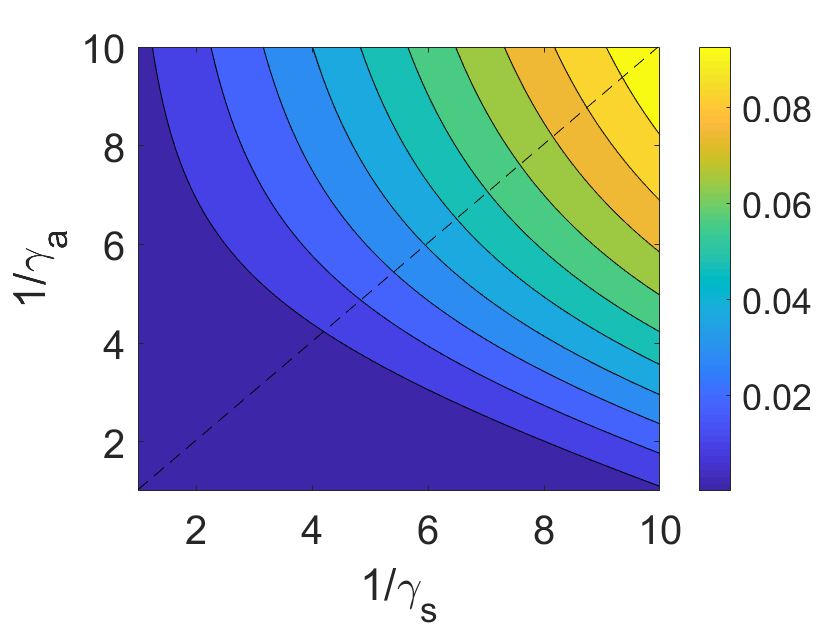}
\caption{}
\end{center}
\end{subfigure}%
~
\begin{subfigure}[b]{0.5\textwidth}
\begin{center}
\includegraphics[width=7cm]{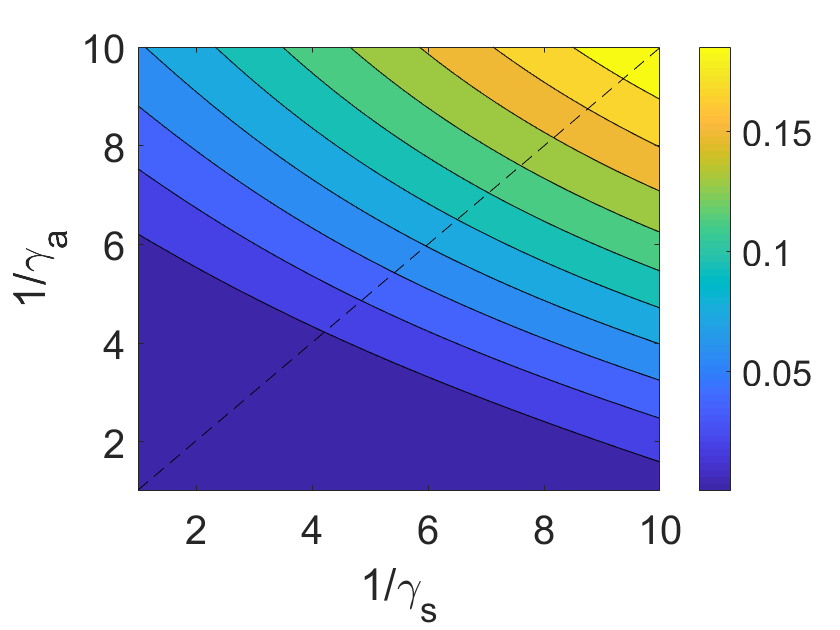}
\caption{}
\end{center}
\end{subfigure}
\caption{\label{different_gamma} Height of the peak of (a,c)  $I_{sQ}$ and (b,d)  $I_a+I_{aQ}+I_{sQ}$ predicted by model \eqref{3}
as a function  of the parameters $1/\gamma_s$ and $1/\gamma_a$. The other parameters are $\psi_s=\psi_a=0$, $\beta_s=\beta_a=0.5$, $S_{sQ}=S_{aQ}=0.1$; the latent period equals   $ \lambda_s=\lambda_a=0.5$  
on panels  (a, b) and 
$\lambda_s=\lambda_a=10$ 
on panels (c,d).
}
\end{figure}

The peak height of the total infected population $I_a + I_{aQ} + I_{sQ}$ is more sensitive to the parameter $1/\gamma_a$ than to the parameter $1/\gamma_s$.  As can be expected, the difference of the sensitivities with respect to these parameters decreases with the increasing latency period when both groups present no symptoms,
see Figures \ref{different_gamma} (b,d). On the other hand, the peak height of $I_{sQ}$ is more sensitive to the parameter $1/\gamma_s$ than to the parameter $1/\gamma_a$, see Figures \ref{different_gamma} (a,c).
This can be attributed to the fact that the number of active symptomatic cases is directly proportional to the average recovery period $1/\gamma_s$, while the undetected asymptomatic cases are not accounted for by $I_{sQ}$.
 
\begin{figure}[h!]
\begin{subfigure}[b]{0.5\textwidth}
\begin{center}
\includegraphics[width=7cm]{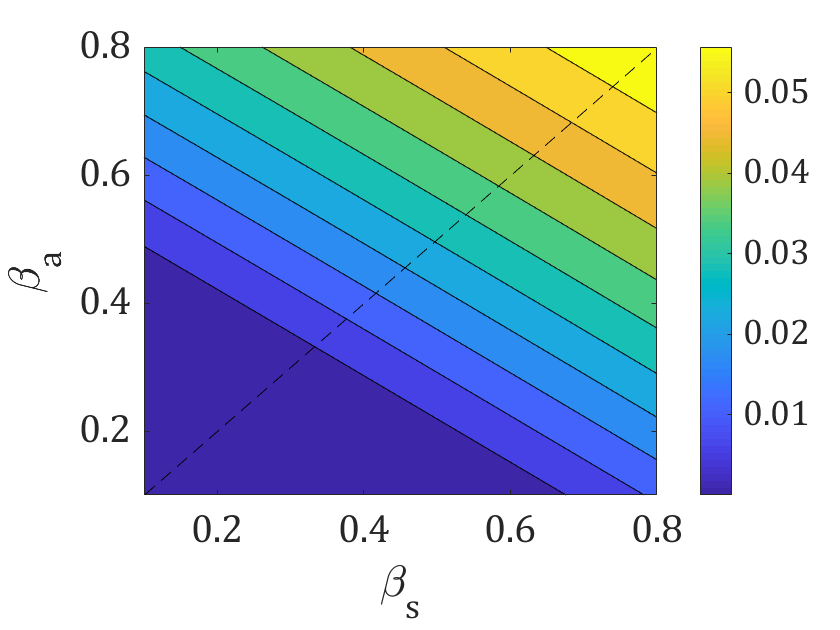}
\caption{}
\end{center}
\end{subfigure}%
~
\begin{subfigure}[b]{0.5\textwidth}
\begin{center}
\includegraphics[width=7cm]{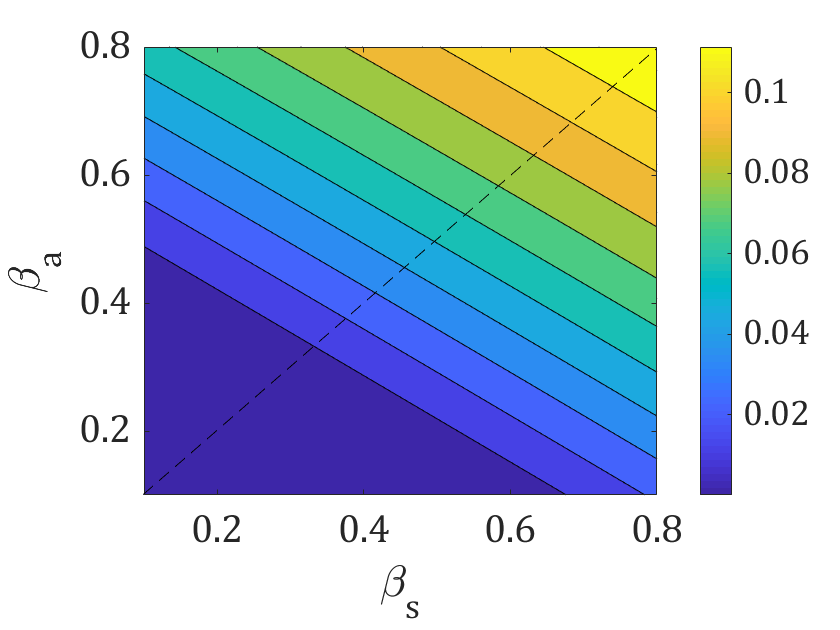}
\caption{}
\end{center}
\end{subfigure}

\begin{subfigure}[b]{0.5\textwidth}
\begin{center}
\includegraphics[width=7cm]{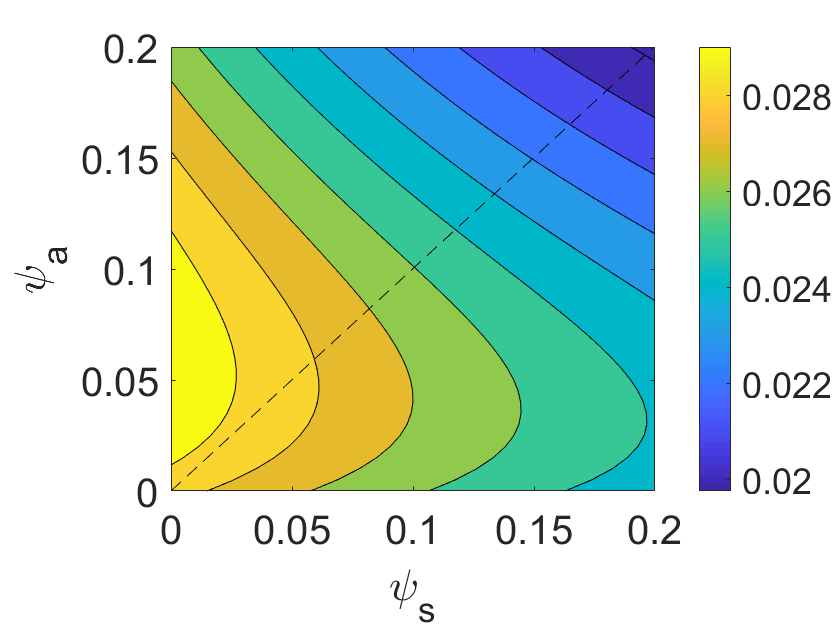}
\caption{}
\end{center}
\end{subfigure}%
~
\begin{subfigure}[b]{0.5\textwidth}
\begin{center}
\includegraphics[width=7cm]{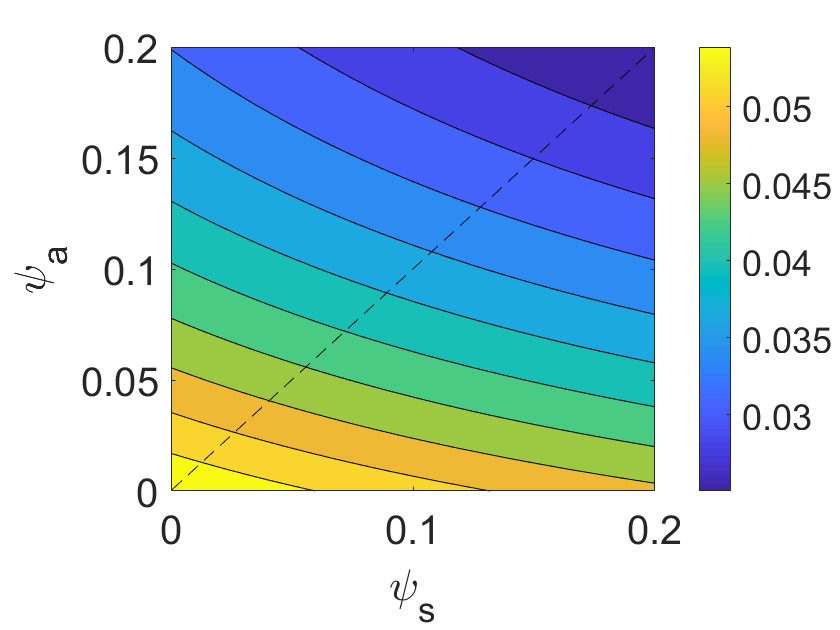}
\caption{}
\end{center}
\end{subfigure}
\caption{\label{different_params}  Height of the peak of (a) $I_{sQ}$ and (b) $I_a+I_{aQ}+I_{sQ}$ 
as a function of $\beta_s$ and $\beta_a$ for $\gamma_s=\gamma_a=1/3.5$, $\psi_s=\psi_a=0$  and $S_{sQ}=S_{aQ}=0.1$. Height of the peak of (c) $I_{sQ}$ and (d) $I_a+I_{aQ}+I_{sQ}$, predicted by model \eqref{3}, as a function of $\psi_s$ and $\psi_a$ for $\gamma_s=\gamma_a=1/3.5$, $\beta_s=\beta_a=0.5$ and $S_{sQ}=S_{aQ}=0.1$. 
}
\end{figure}

It is a controversial topic whether asymptomatic individuals spread the virus less than symptomatic individuals do. Figure \ref{different_params} (a,b) shows that
 the peaks of $I_{sQ}$ and $I_a + I_{aQ} + I_{sQ}$ have higher sensitivity to the parameter $\beta_a$ than to $\beta_s$
 when the other parameters of the symptomatic and asymptomatic groups are approximately equal.
 

\begin{figure}[h!]
\begin{center}
\includegraphics[width=10cm]{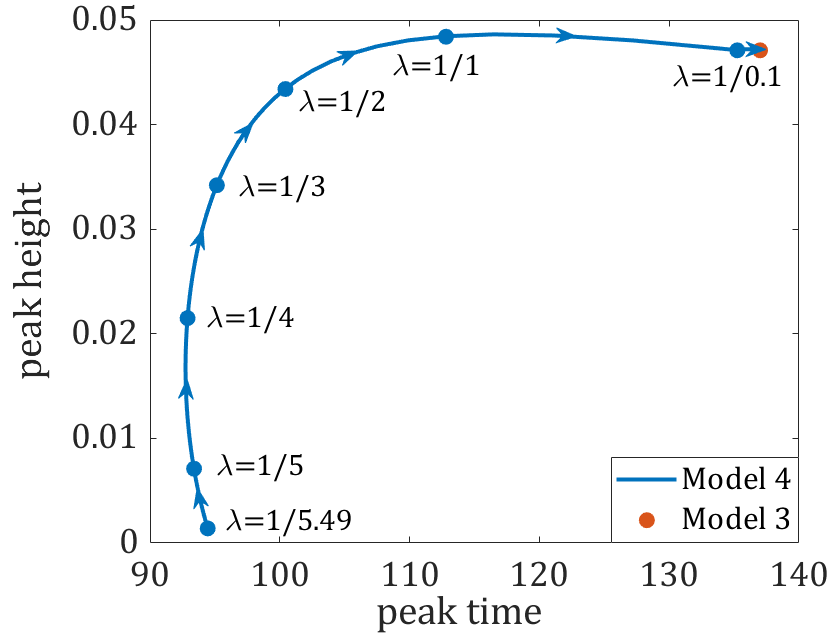}
\caption{\label{comparison} Dependence of timing and height of the peak of $I_{sQ}$ on $\lambda_s$ in model \eqref{3} (blue points)   for $1/\lambda_s+1/\gamma_s =  1/\delta=5.5$ and $1/\lambda_s \in (0,5.5).$ Here $\lambda_s=\lambda_a$, $\gamma_s=\gamma_a$, $\psi_s=\psi_a=0$, $E_s(0)=E_a(0)=1/N$ and $S_{sQ}=S_{aQ}=0$.
The red point corresponds to model \eqref{model1}.}
\end{center}
\end{figure}

One can assume that asymptomatic individuals have little reason to test unless required to do so.
Further, according to \cite{psi}, asymptomatic individuals can have a weaker immune response to COVID-19 infection,  
hence symptomatic cases (even during the latent phase) might be better detectable
than asymptomatic ones. These considerations suggest that $\psi_a<\psi_s$.
%
Figure \ref{different_params} (c,d) shows  
that in order to ensure that the peak of $I_{sQ}$ and $I_{sQ}+ I_{aQ}+ I_a$ is low (the dark blue region), it is necessary to effectively detect and quarantine not only symptomatic individuals in the latent phase of the disease but also a relatively high proportion of asymptomatic cases. The peak of $I_{sQ}+ I_{aQ}+ I_a$ is more sensitive to the testing rate
$\psi_a$ than to $\psi_s$.
The non-monotonic dependence of the peak of $I_{sQ}$ on $\psi_a$ in  Figure \ref{different_params} (c)
is due to the fact that 
an increase of the testing rate (for small values of the testing rate) results in 
an increase of the number of reported cases due to a
higher detection rate (the same effect is observed in Figure \ref{PSI1} (a,b)).

 The symptomatic cases are detected and quarantined when the symptoms present themselves after the latent phase.
 On the other hand, the asymptomatic  infected  individuals, if undetected by testing, continue spreading the virus 
 after the latent phase during the whole infectious stage. This can explain the higher sensitivity of the peak of the total infected population to the parameters of the asymptomatic group, as Figures \ref{different_gamma}, \ref{different_params} demonstrate for the transmission, recovery and testing rates. Since 
 the estimates of the epidemiological parameters are currently based mainly on the data obtained from symptomatic cases,
 the above sensitivity analysis suggests that better estimates of the parameters of the asymptomatic group
 can play an important role for improving the predictive power of models of the epidemic dynamics.



The duration of the latent phase can have a significant impact on 
dynamics of $I_{sQ}$.
Figure \ref{comparison} shows the timing and height of the  peak of  $I_{sQ}$ for model \eqref{3} for different values of $\lambda=\lambda_s=\lambda_a=\lambda$ under the assumption that $\gamma_s=\gamma_a$ and $1/\lambda_s + 1/\gamma_s = 1/\delta$ with a fixed $\delta$.   
Both the height and timing of the peak increase as the duration $1/\lambda$ of the latent phase decreases.
As the length of the latent phase $1/\lambda$ approaches  zero, the dynamics of model \eqref{3} converges to the dynamics of model \eqref{model1}. 


\subsection{
Indiscriminate quarantining vs massive testing}

\begin{figure}[h!]
\begin{subfigure}[b]{0.5\textwidth}
\begin{center}
\includegraphics[width=7cm]{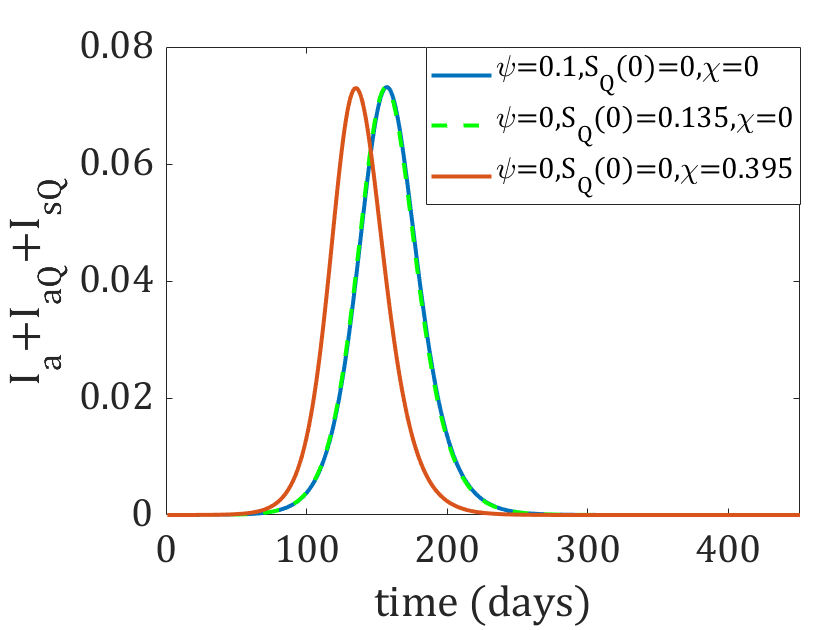}
\caption{}
\end{center}
\end{subfigure}%
~
\begin{subfigure}[b]{0.5\textwidth}
\begin{center}
\includegraphics[width=7cm]
{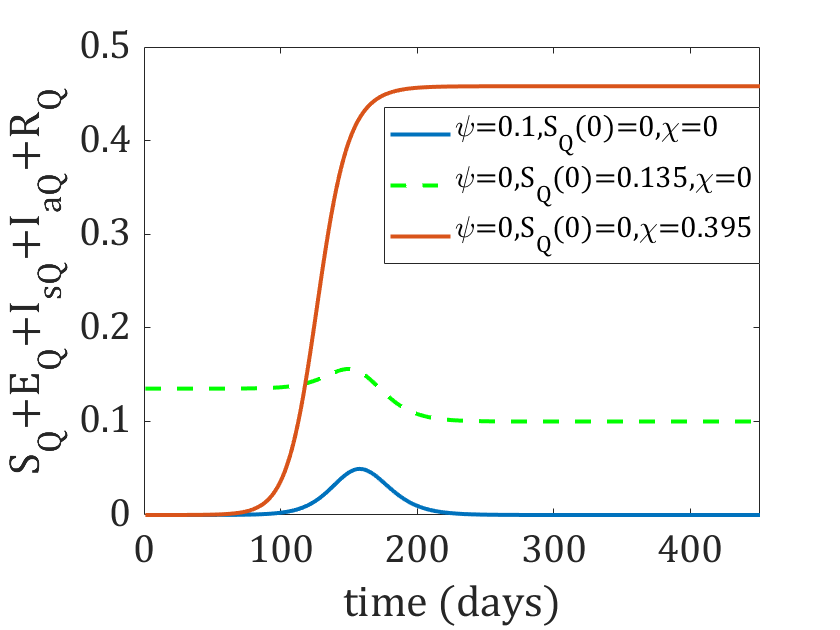}
\caption{}
\end{center}
\end{subfigure}

\begin{subfigure}[b]{0.5\textwidth}
\begin{center}
\includegraphics[width=7cm]{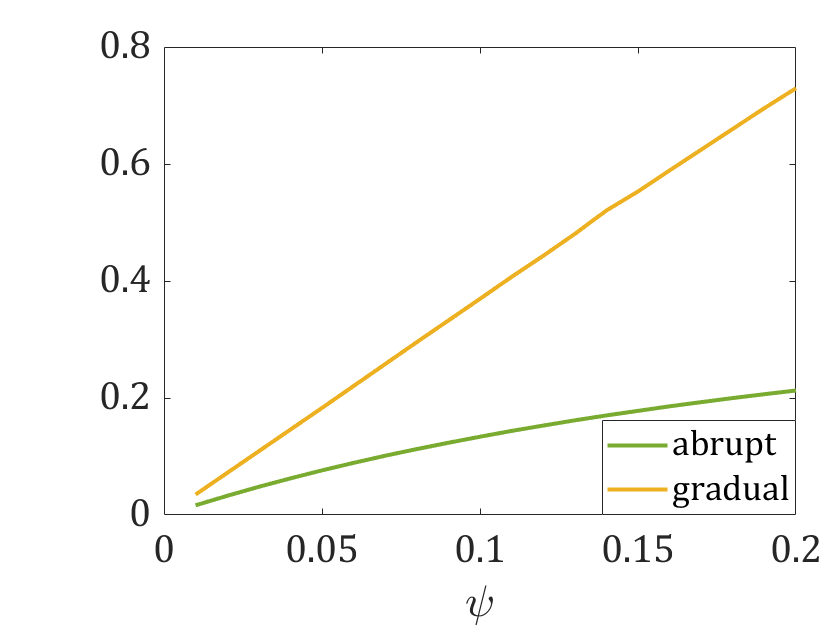}
\caption{}
\end{center}
\end{subfigure}%
~
\begin{subfigure}[b]{0.5\textwidth}
\begin{center}
\includegraphics[width=7cm]{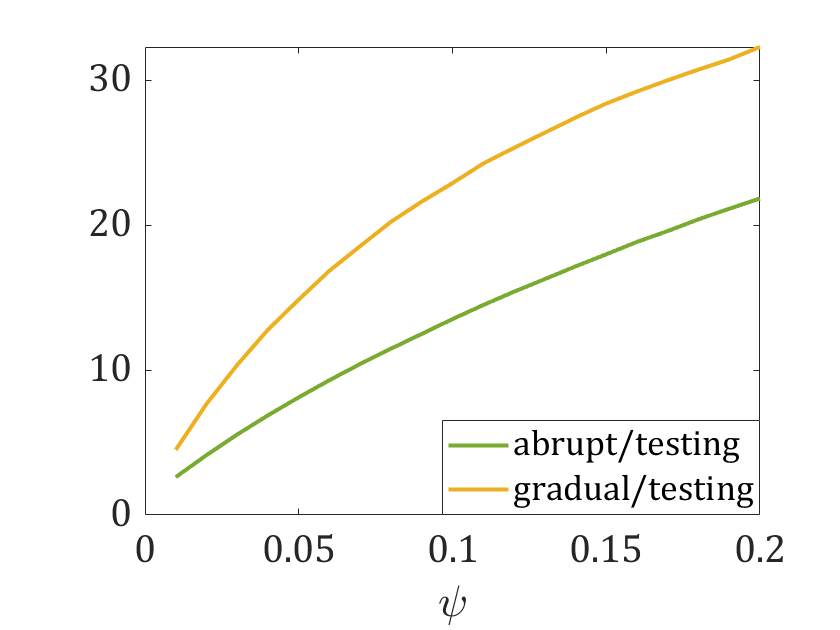}
\caption{}
\end{center}
\end{subfigure}

\caption{Indiscriminate quarantining versus quarantining based on testing in model \eqref{model1}  with $E_Q(0)=2/N.$ (a) Dynamics of the total infected population,  $I_a+I_{aQ}+I_{sQ}$,  for a selected testing rate $\psi$ and the corresponding value of $S_Q(0)$ and $\chi$ resulting in the same height of the peak of $I_a+I_{aQ}+I_{sQ}$. (b) Dynamics of the total quarantined population corresponding to the graphs on panel (a); the same color coding is used.
(c) For every point $(\psi_0, s_0)$ of the green curve, the testing rate $\psi=\psi_0$ with zero indiscriminate quarantining 
($S_Q(0)=0$, $\chi=0$) results 
in the same peak of infection $I_a+I_{aQ}+I_{sQ}$ 
as the initial quarantine $S_Q(0)=s_0$ with zero testing ($\psi=\chi=0$).
For every point $(\psi_0, \chi_0)$ of the orange curve, the testing rate $\psi=\psi_0$ with zero indiscriminate quarantining 
($S_Q(0)=0$, $\chi=0$) results 
in the same peak of infection 
as the indiscriminate quarantining rate $\chi=\chi_0$ with zero testing and zero initial quarantine ($\psi=0,$ $S_Q(0)=0$).
(d) The ratio of the time-integral of the total quarantined population in the case of indiscriminate quarantining
to the value of this integral in the case of quarantining based on testing under the constraint that the infection peak is the same for both scenarios. 
The green and orange curves show the dependence of this ratio on $\psi$ and
correspond to the curves of the same colors, respectively, on panel (c).
\label{PSI_SQ}}
\end{figure}

Due to economic and social reasons, it is clear that  lockdowns  applied in many countries cannot last for very long.  Therefore there is a  demand to find a more economically sustainable solution to replace massive home quarantine measures.  One reasonable possibility seems to be  effective testing.  
The blue curve in Figure \ref{PSI_SQ} (a) shows the evolution of the infected population for model \eqref{model1} with 13.5\% of the total population quarantined at the initial moment and zero testing rate;  the green curve corresponds to the situation where nobody stays in quarantine initially  but 
10\% 
of the total population is screened daily for COVID-19 (more precisely, 10\% of infected asymptomatic individuals are successfully detected daily).
Under both  setups the total infected population $I_a+I_{aQ}+I_{sQ}$ follows approximately the same trajectory. 
However, comparing the dynamics of the quarantined population, we see that
it is significantly smaller at all times when testing is used than 
in the case when indiscriminate quarantining is applied, see the green and blue curves in Figure \ref{PSI_SQ} (b).
A similar observation can be made when instead of abrupt quarantining the gradual quarantining strategy
is applied, see the red curves in Figure \ref{PSI_SQ} (a,b). Again, the quarantined population is significantly smaller
when quarantining is based on testing than in the case when indiscriminate gradual quarantining is used, while the trajectories of the infected population are similar for both scenarios.
Recall that the number of active symptomatic cases is proportional to the total infected population and equals $(1-k)(I_a+I_{aQ}+I_{sQ})$.

Figure \ref{PSI_SQ} (c,d) presents similar results 
for  a range of values  of the testing rate. We compare the scenario with a positive testing rate $\psi$,
zero indiscriminate quarantining rate $\chi=0$ and zero initial quarantine size $S_Q(0)=0$ to the scenarios
with $\psi=\chi=0$, $S_Q(0)>0$ and with $\psi=0$, $\chi>0$, $S_Q(0)=0$.
The first scenario corresponds to quarantining based on testing, while the
second and third scenarios correspond to indiscriminate abrupt and gradual quarantining, respectively.
Given a $\psi>0$ in the first scenario, the $S_Q(0)>0$ and $\chi>0$ in the second and third scenarios are selected in such a way as to ensure that the height of the infection peak in all the three scenarios is the same (as in Figure \ref{PSI_SQ} (a)).
This constraint defines $S_Q(0)$ and $\chi$ in the second and third scenarios, respectively, as functions of the testing rate $\psi$ used in the first scenario, see the graphs in Figure \ref{PSI_SQ} (c).
For each point of these graphs, we also compute the time-integral of the total quarantined population over the duration of the epidemic; the latter is technically defined as the the time interval $[0,t_e]$ where $t_e$ is the moment after the infection peak when $I_{sQ}(t_e)=0.1/N$. Figure \ref{PSI_SQ} (d) shows that the time-integral of the quarantined population is by an order of magnitude smaller when quarantining is based on massive testing than in the case when indiscriminate (abrupt or gradual) quarantining is applied. Moreover, the ratio of the value of this time-integral in the case of indiscriminate quarantining
to its value in the case of quarantining based on testing increases with the increasing testing rate.

\begin{figure}[h]
\begin{center}
\begin{subfigure}[b]{0.5\textwidth}
\begin{center}
\includegraphics[width=7.2cm]{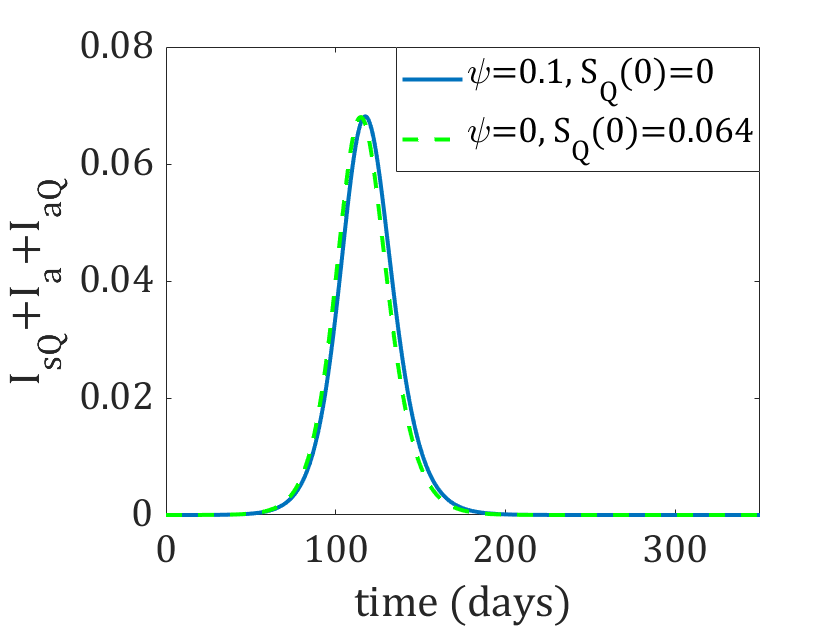}
\caption{}
\end{center}
\end{subfigure}%
~
\begin{subfigure}[b]{0.5\textwidth}
\begin{center}
\includegraphics[width=7.2cm]{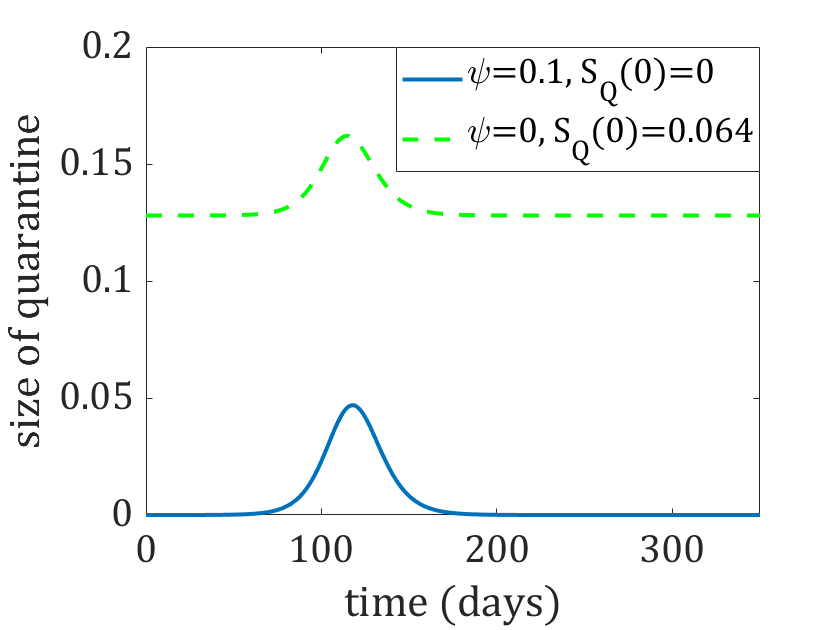}
\caption{}
\end{center}
\end{subfigure}
\caption{\label{psi_SQ2}Indiscriminate quarantining versus quarantining based on testing in model \eqref{3}  with $S_{sQ}(0)=S_{aQ}(0)=S_Q(0)$ and $E_s(0)=E_a(0)=1/N$ and $\psi_a=\psi_s=\psi$. (a) Dynamics of the total infected population,  $I_a+I_{aQ}+I_{sQ}$,  for a selected testing rate $\psi$ and the corresponding value of $S_Q(0)$  resulting in the same height of the peak of $I_a+I_{aQ}+I_{sQ}$. (b) Dynamics of the total quarantined population corresponding to the graphs on panel (a); the same color coding is used.
}
\end{center}
\end{figure}

Similar results were obtained with model \eqref{3}, see Figure \ref{psi_SQ2}.

These observations indicate that effective testing might be a way to replace massive quarantine measures. Moreover, they may help to explain why  some countries such as South Korea and Singapore, which implemented 
massive rapid free testing and an extremely good system of contact tracing,
have been most successful in containing the
COVID-19 outbreak.

\section{Conclusions}
Several issues involved in  modeling  the spread of an epidemic in a population with symptomatic and asymptomatic groups were discussed. We introduced two models which differ in their complexity and assumptions.  Both models account for two intervention strategies\,---\,indiscriminate home quarantining of part of the population and isolation of positively tested individuals.
In the  numerical analysis we used average values of epidemiological parameters according to the current publications on COVID-19.

The models predict that dynamics of infection are similar for both intervention strategies.
However, massive testing allows to control and contain the infection spread using a much lower isolation rate than in the case of indiscriminate quarantining. In particular, given a constraint that limits the maximal number of simultaneous active cases,
the isolation rate
decreases with the increasing testing rate.
The models predict that massive testing measures can help to decrease
the size of the quarantine, which enforces the above constraint, 
by an order of magnitude and more. 
Due to a variety of complicating factors,
these findings cannot be considered as a quantitatively accurate prediction for the coronavirus pandemic. However, they 
indicate that massive testing can potentially stand-in for costly quarantine measures.
As such, our findings agree with the earlier predictions based on different models 
\cite{testing}.
 


We also discussed the role of asymptomatic individuals for the spread of the epidemic. 
Our simulations  indicate that a possible difference of the epidemiological parameters of the symptomatic and asymptomatic groups
can play an important role for predicting
 the infection peak value and duration of the epidemic. 
We showed that the peak of the infected population is more sensitive to the 
rate of transmission from the asymptomatic group
and to the rate of testing of this group (for higher rates of testing) than to the corresponding parameters of the symptomatic group.
This can be attributed to the fact that symptomatic cases are  detected and quarantined 
earlier, after the latent phase of the disease, while
the asymptomatic  infected  individuals, if undetected by testing, continue spreading the virus 
during the whole infectious phase.  
Therefore  more accurate estimates of the epidemiological  parameters of COVID-19, especially those of the asymptomatic groups, can improve predictions of the epidemic dynamics.




\end{document}